\documentclass[12pt]{article}
\usepackage {amsmath, amssymb}
\def\ep{\varepsilon}
\def\th{\vartheta}
\begin{document}

\title{{\large
\bf Exact $L_2$-small ball asymptotics of Gaussian processes \\ and the
spectrum of boundary value problems \\ with "non-separated"\ boundary
conditions
}}
\author { {\it Alexander I.~Nazarov}\footnote
{Partially supported by RFFR grant No.07-01-00159}
\\ Dept. of Mathematics and Mechanics,\\
St.Petersburg State University, 198504, Russia\\
{\small e-mail:\ an@AN4751.spb.edu}} \maketitle
\begin{abstract}

{\it\small We sharpen a classical result on the spectral asymptotics of the
boundary value problems for self-adjoint ordinary differential operator.
Using this result we obtain the exact $L_2$-small ball asymptotics for a new
class of zero mean Gaussian processes. This class includes, in particular,
integrated generalized Slepian process, integrated centered Wiener process
and integrated centered Brownian bridge.}
\end{abstract}

\section*{Introduction}

\indent The problem of small ball behavior for norms of Gaussian processes
was actively studied in recent years, see, for example, the reviews \cite{Lf}
and \cite{LS}. We discuss the most explored case of $L_2$-norm. Suppose we
have a Gaussian process $X(t),0\le t\le 1$, with zero mean and covariance
function $G_X(t,s)=EX(t)X(s)$, $t,s\in [0,1]$. Let $\|X\|=\|X\|_{L_2(0,1)}$
and consider
$$
Q(X;\ep)={\bf P}\{\|X\|\le \ep \}.
$$
The problem to define the behavior of $Q(X;\ep)$ as $\ep\to 0$ was solved in
\cite{S}, but in an implicit way. Therefore, a number of papers provided the
simplification of the expression for $Q(X;\ep )$ under various assumptions
(see, e.g., the references in \cite{LS} and in \cite{DLL}).

According to the classical Karhunen- Lo\`eve expansion one has for the
process $X$ the equality in distribution
$$\|X\|^2 =\int\limits_0^1X^2(t)\,dt\stackrel {d}{=}
\sum_{n=1}^\infty \lambda_n\eta_n^2,$$

\noindent where  $\eta_n$, $n\in\mathbb N$, are independent standard
Gaussian random variables while $\lambda_n=\lambda_n(X)>0$, $n\in\mathbb N$,
$\sum\limits_n\lambda_n <\infty $ are the eigenvalues of the integral
equation
$$ \lambda y(t)=\int\limits_0^1
G_X(t,s)y(s)\,ds,\quad 0\le t\le1.$$

Thus we are led to the equivalent problem of studying the asymptotics of
${\bf P}\left\{\sum_{n=1}^\infty \lambda _n \eta_n^2\le\ep ^2 \right\}$ as
$\ep \to 0$. Unfortunately, explicit formulas for eigenvalues can be obtained
only in a limited number of examples.

A new approach developed in the paper \cite{NN} gives exact (up to a
constant) small ball asymptotics for Gaussian process $X$ under assumption
that $G_X$ is the Green function of a boundary value problem (BVP) for
ordinary differential operator with "separated"\ boundary conditions
(Sturm-type conditions). This work was completed by the paper \cite{Na} where
sharp constants in the small ball asymptotics were calculated for many
Gaussian processes. A part of results of \cite{Na} was independently obtained
in \cite{GHLT1}, \cite{GHLT2}. In \cite{NP} these results were transferred to
a class of weighted processes.

The approach of \cite{NN} is based on classical Birkhoff's results on
the spectral asymptotics of BVPs to ordinary differential operators. It is
well-known, see, e.g., \cite[\S4]{Nm}, \cite[Ch.XIX]{DS}, that eigenvalues
$\mu_n$ of regular (in particular, self-adjoint) BVPs can be expanded into
asymptotic series in powers of $n$. The first term of this expansion is
completely determined by the main coefficient of the operator while the
formulas for other terms are rather complicated. It turned out that in the
case of "separated"\ boundary conditions the second term of the asymptotics
is completely determined by the sum of orders of boundary functionals.
Therefore, it can be derived in explicit form without additional assumptions.
Having in hands two-term asymptotics for $\mu_n$ (and consequently for
$\lambda_n=\mu_n^{-1}$) we can apply the approach from \cite{DLL} and
comparison theorem \cite{Li} to obtain the final result.\medskip

In general case, where the boundary conditions are non-separated, the
eigenvalues of BVP can be split into two subsequences, and the formulas
for the second terms of the asymptotics of these subsequences, generally
speaking, cannot be simplified. However, the Lifshits lemma (see below),
combined with \cite[Theorem 6.2]{NN}, shows that the $L_2$-small ball
behavior up to constant for the corresponding Gaussian process depends
only on {\bf the sum} of these second terms. In this paper we show that this
sum, as before, is completely determined by the sum of orders of boundary
functionals. This allows us to generalize the results of \cite{NN} to
considerably larger class of processes.

As for sharp constants in the small ball asymptotics, one can write down
explicit formulas for them if the eigenfunctions of the covariance kernel
can be expressed in terms of elementary or special functions. In this case
the asymptotics of corresponding Fredholm determinants can be calculated by
the complex variable methods, see \cite{Na}, \cite{GHLT1}. In this paper we
show this by example of several well-known processes generating boundary
value problems with non-separated boundary conditions.\medskip

The paper is organized as follows. In Section 1 we prove the theorem on the
second terms of spectral asymptotics to BVPs with non-separated boundary
conditions. Also the small ball asymptotics up to constant for corresponding
Gaussian processes is given. The sharp small deviation constants for multiply
integrated generalized Slepian process, for some variants of integrated
centered Brownian bridge and for some kinds of integrated centered Wiener
process are calculated, respectively, in Sections 2, 3 and 4. We note that
the $L_2$-small ball asymptotics for some centered processes was derived in
\cite{BNO}.\medskip

Let us recall some notation. For any zero mean Gaussian process $X(t)$,
$0\le t\le 1$, we introduce {\bf the centered process}
${\overline X}(t)=X(t)-\int_0^1X(s) \,ds$ and {\bf the $m$-times integrated}
process
$$X_m^{[\beta_1,\,...,\,\beta_m]}(t) =(-1)^{\beta_1+\,\dots\,+\beta_m}
\underbrace{\int
_{\beta_m}^t\dots\int
_{\beta_1}^{t_1}}_{m}X(s)\,ds\,dt_1\dots$$
(here any index $\beta_j$ equals either zero or one, $0\le t\le1$). For the
sake of brevity the upper index is sometimes omitted.

The function $G(t,s)$ is called the Green function of (self-adjoint) boundary
value problem for differential operator $L$ if it satisfies the equation
$LG=\delta(t-s)$ in the sense of distributions and satisfies the boundary
conditions. The existence of Green function is equivalent to the
invertibility of operator $L$ with given boundary conditions, and $G(t,s)$ is
the kernel of the integral operator $L^{-1}$. If homogeneous BVP has a
non-trivial solution $\varphi_0$ (without loss of generality it can be
assumed to be normalized in $L_2(0,1)$) then the Green function obviously
does not exist. If $\varphi_0$ is the unique solution up to a constant
multiplier\footnote {In a similar way one can consider the case of the
multiple zero eigenvalue but we do not need it.} then the function $G(t,s)$
is called {\bf generalized Green function} provided it satisfies the equation
$LG=\delta(t-s)-\varphi_0(t)\varphi_0(s)$ in the sense of distributions,
satisfies the boundary conditions and the orthogonality condition
$$\int\limits_0^1G(t,s)\varphi_0(s)\,ds=0 \quad \mbox {for all} \quad
0\le t\le 1.\eqno(0.1)$$
The generalized Green function is the kernel of the integral operator which
is inverse to $L$ on the subspace of functions orthogonal to $\varphi_0$ in
$L_2(0,1)$. The reader is referred to \cite[Chapter 2, \S1]{Sm} for more
detailed properties of the Green function and the generalized Green function
(for the second order operators).

The space $W_p^m(0,1)$ is the Banach space of functions $u$ having
continuous derivatives up to $(m-1)$-th order when $u^{(m-1)}$ is
absolutely continuous on $[0,1]$ and $u^{(m)}\in L_p(0,1)$.
If $p=2$ it is a Hilbert space.

We set $z_{\ell}=\exp(i\pi/\ell)$ while $\mathfrak V(...)$ stands for the
Vandermond determinant:
$$
\mathfrak V(\alpha_1,\alpha_2,\ldots,\alpha_n)=\det
\begin{bmatrix}
1 & \alpha_1 & \alpha_1^2 & \ldots & \alpha_1^{n-1}\\
1 & \alpha_2 & \alpha_2^2 & \ldots & \alpha_2^{n-1}\\
\hdotsfor{5}\\
1 & \alpha_n & \alpha_n^2 & \ldots & \alpha_n^{n-1}&
\end{bmatrix}
=\prod\limits_{1\le j<k\le n}\!(\alpha_k-\alpha_j).
$$

We cite the statement due to M.A. Lifshits, see \cite{BNO}. This statement is
repeatedly used in our paper.\medskip

{\bf Lemma 0.1.} { \it Let $V_1,V_2>0$ be two independent random variables
with given small ball behavior; namely, let as $r\to 0$
$$
{\bf P}\{V_1\le r\}\sim K_1r^{a_1}\exp (-D_1^{d+1}r^{-d}),\qquad
{\bf P}\{V_2\le r\}\sim K_2r^{a_2}\exp (-D_2^{d+1}r^{-d}).
$$
Then their sum has the following small ball asymptotics:
$$
{\bf P}\{V_1+V_2\le r\}\sim Kr^a\exp (-D^{d+1}r^{-d}),
$$
where}
$$
D=D_1+D_2,\qquad
a=a_1+a_2-{\textstyle\frac d2},\qquad
K=K_1K_2\sqrt {\frac {2\pi d}{d+1}}\cdot
\frac {D_1^{a_1+\frac 12}D_2^{a_2+\frac 12}}{D^{a+\frac 12}}. $$

\section{Eigenvalues asymptotics for BVPs and\\ small ball asymptotics}

Let $\cal L$ be a self-adjoint differential operator of order $2\ell$
generated by a differential expression
$${\cal L}u\equiv
(-1)^{\ell}\left(p_{\ell}u^{({\ell})}\right)^{({\ell})}+
\left(p_{\ell-1}u^{({\ell}-1)}\right)^{({\ell}-1)}+\dots+p_0u,\eqno(1.1)$$
($p_{\ell}(x)>0$) and by $2\ell$ boundary conditions
$$U_{\nu}(u)\equiv U_{\nu 0}(u)+U_{\nu 1}(u)=0,\qquad \nu=1,\dots,2\ell,
\eqno(1.2)$$
where
$$U_{\nu 0}(u)=\alpha_{\nu}u^{(k_{\nu})}(0)+\sum\limits_{j=0}^{k_{\nu}-1}
\alpha_{\nu j}u^{(j)}(0),$$
$$U_{\nu 1}(u)=\gamma_{\nu}u^{(k_{\nu})}(1)+\sum\limits_{j=0}^{k_{\nu}-1}
\gamma_{\nu j}u^{(j)}(1),$$
and for any index $\nu$ at least one of coefficients $\alpha_{\nu}$ and
$\gamma_{\nu}$ is not equal to zero.

It is well known, see, e.g., \cite[\S4]{Nm}, that the system of boundary
conditions (1.2) can be reduced to the {\bf normalized form} by equivalent
transformations. In what follows we always assume that this reduction is
realized. This form is specified by the minimal sum of orders of all boundary
conditions. Since this quantity is of great importance in our arguments, we
introduce the notation $\varkappa=\sum\limits_{\nu=1}^{2\ell} k_{\nu}$. We
remark also that the inequalities
$$2\ell-1\ge k_1\ge k_2\ge\dots\ge k_{2\ell}\ge0,\qquad k_{\nu}>k_{\nu+2}$$
hold true for the normalized boundary conditions.

For simplicity we assume $p_j\in W^j_{\infty}[0,1]$, $j=0,\dots,\ell$.
Then the domain ${\cal D}({\cal A})$ consists of the functions
$u\in W^{2\ell}_2(0,1)$ satisfying boundary conditions (1.2).\medskip

Consider the eigenvalue problem
$${\cal L}u=\mu u\quad \mbox{\rm on} \quad [0,1], \qquad
u\in {\cal D} ({\cal L}).\eqno(1.3)$$
It is well known, see, e.g., \cite[\S4, Theorem 2]{Nm}, that for
$p_{\ell}\equiv1$ the eigenvalues of (1.3) counted according to their
multiplicities can be split into two subsequences $\mu'_n$, $\mu''_n$,
$n\in{\mathbb N}$, such that, as $n\to\infty$,
$$\mu'_n=\left(2\pi n+\rho'+O(n^{-1/2})\right)^{2\ell}, \qquad
\mu''_n=\left(2\pi n+\rho''+O(n^{-1/2})\right)^{2\ell},\eqno(1.4)$$
where $\xi'=\exp(i\rho')$, $\xi''=\exp(i\rho'')$ are the roots of quadratic
equation
\begin{multline*}
\theta_1\xi+\theta_0+ \theta_{-1}\xi^{-1}\equiv\\
\equiv\det
\begin{bmatrix}
(\alpha_1+\xi \gamma_1)&\alpha_1\omega_1^{k_1}&\dots&
\alpha_1\omega_{\ell-1}^{k_1}&\omega_{\ell}^{k_1}(\alpha_1+\xi^{-1}\gamma_1)&
\gamma_1\omega_{\ell+1}^{k_1} &\dots&\gamma_1\omega_{2\ell-1}^{k_1}\\
(\alpha_2+\xi \gamma_2)&\alpha_2\omega_1^{k_2}&\dots&
\alpha_2\omega_{\ell-1}^{k_2}&\omega_{\ell}^{k_2}(\alpha_2+\xi^{-1}\gamma_2)&
\gamma_2\omega_{\ell+1}^{k_2} &\dots&\gamma_2\omega_{2\ell-1}^{k_2}\\
\hdotsfor[3]{8} \\
\hdotsfor[3]{8} \\
(\alpha_{2\ell}+\xi \gamma_{2\ell})&\alpha_1\omega_1^{k_{2\ell}}&\dots&
\alpha_{2\ell}\omega_{\ell-1}^{k_{2\ell}}&
\omega_{\ell}^{k_{2\ell}}(\alpha_{2\ell}+\xi^{-1}\gamma_{2\ell})&
\gamma_{2\ell}\omega_{\ell+1}^{k_{2\ell}} &\dots&
\gamma_{2\ell}\omega_{2\ell-1}^{k_{2\ell}}
\end{bmatrix}=0
\end{multline*}

\noindent (here $\omega_j=z_{\ell}^j=\exp(ij\pi/\ell)$).\medskip

{\bf Remark 1}. In general case the problem (1.3) can be reduced to the case 
$p_{\ell}\equiv1$ by the independent variable transform, see \cite[\S4]{Nm}.
The expressions in brackets in (1.4) should be divided by
$\int\limits_0^1 p_{\ell}^{-1/(2\ell)}(x)\,dx$.\medskip

{\bf Theorem 1.1}.{ \it The following relation holds true:}
$$\rho'+\rho''=2\pi\ell-3\pi-\frac {\pi\varkappa}{\ell}.\eqno(1.5)$$

\noindent {\bf Proof}. The key observation is that
$\xi\xi'=\frac{\theta_{-1}}{\theta_1}$ does not depend on $\alpha_{\nu}$ and
$\beta_{\nu}$ as long as $\theta_1\ne0$ i.e. as long as boundary conditions
(1.2) are regular. To check this fact we write
$$\theta_1=\det\begin{bmatrix}
\gamma_1&\alpha_1\omega_1^{k_1}&\dots&\alpha_1\omega_{\ell-1}^{k_1}&
\alpha_1\omega_{\ell}^{k_1}&\gamma_1\omega_{\ell+1}^{k_1} &\dots&
\gamma_1\omega_{2\ell-1}^{k_1}\\
\gamma_2&\alpha_2\omega_1^{k_2}&\dots&\alpha_2\omega_{\ell-1}^{k_2}&
\alpha_2\omega_{\ell}^{k_2}&\gamma_2\omega_{\ell+1}^{k_2} &\dots&
\gamma_2\omega_{2\ell-1}^{k_2}\\
\hdotsfor[3]{8} \\
\hdotsfor[3]{8} \\
\gamma_{2\ell}&\alpha_1\omega_1^{k_{2\ell}}&\dots&
\alpha_{2\ell}\omega_{\ell-1}^{k_{2\ell}}&
\alpha_{2\ell}\omega_{\ell}^{k_{2\ell}}&
\gamma_{2\ell}\omega_{\ell+1}^{k_{2\ell}} &
\dots&\gamma_{2\ell}\omega_{2\ell-1}^{k_{2\ell}}
\end{bmatrix},
$$
$$\theta_{-1}=\det\begin{bmatrix}
\alpha_1&\alpha_1\omega_1^{k_1}&\dots&\alpha_1\omega_{\ell-1}^{k_1}&
\gamma_1\omega_{\ell}^{k_1}&\gamma_1\omega_{\ell+1}^{k_1} &\dots&
\gamma_1\omega_{2\ell-1}^{k_1}\\
\alpha_2&\alpha_2\omega_1^{k_2}&\dots&\alpha_2\omega_{\ell-1}^{k_2}&
\gamma_2\omega_{\ell}^{k_2}&\gamma_2\omega_{\ell+1}^{k_2} &\dots&
\gamma_2\omega_{2\ell-1}^{k_2}\\
\hdotsfor[3]{8} \\
\hdotsfor[3]{8} \\
\alpha_{2\ell}&\alpha_1\omega_1^{k_{2\ell}}&\dots&
\alpha_{2\ell}\omega_{\ell-1}^{k_{2\ell}}&
\gamma_{2\ell}\omega_{\ell}^{k_{2\ell}}&
\gamma_{2\ell}\omega_{\ell+1}^{k_{2\ell}} &
\dots&\gamma_{2\ell}\omega_{2\ell-1}^{k_{2\ell}}
\end{bmatrix}.
$$
Taking the common multiplier $\omega_1^{k_j}$ over from the $j$-th row in the
first determinant we obtain $\theta_1=-\omega_1^{\varkappa}\theta_{-1}$ and,
therefore,
$$\rho+\rho'=-\frac {\pi\varkappa}{\ell}+(2N-1)\pi,\qquad N\in {\mathbb Z}.
\eqno(1.6)$$

Since the eigenvalues of the problem (1.3) depend on $\alpha_{\nu}$ and
$\gamma_{\nu}$ continuously, the parameter $N$ in (1.6) does not depend on
$\alpha_{\nu}$ and $\gamma_{\nu}$. However, Theorem 7.1 \cite{NN} shows that
if the boundary conditions are {\bf separated}, i.e. $U_{\nu1}\equiv0$
for $\nu=1,\dots,\ell$ and $U_{\nu0}\equiv0$ for $\nu=\ell+1,\dots,2\ell$
then in (1.4)
$$\rho'=\pi(\ell-1)-\frac {\pi\varkappa}{2\ell},\qquad
\rho''=\pi(\ell-2)-\frac {\pi\varkappa}{2\ell}$$
(we note that in this case two sequences $\mu'_n$ and $\mu''_n$ can be
naturally merged into $\mu_n=\left(\pi (n+\ell-1-\frac {\varkappa}{2\ell})
+O(n^{-1/2})\right)^{2\ell}$). Thus, in this case (1.5) holds true, and the
statement follows. \hfill$\square$\medskip

{\bf Remark 2}. In fact the assumption of self-adjointness of the operator
can be relaxed to the assumption of regularity of the system (1.2). The
requirements on the coefficients $p_j$ also can be weakened.\medskip

The next theorem generalizes \cite[Theorem 7.2]{NN} where this result was
obtained for the case of separated boundary conditions.\medskip

{\bf Theorem 1.2}.{ \it Let the covariance $G_X(t,s)$ of a zero mean
Gaussian process $X(t)$, $0\le t \le 1$, be the Green function of a
self-adjoint positive definite operator ${\cal L}_X$ generated by a
differential expression (1.1) and by boundary conditions (1.2). Let
$\varkappa<2\ell^2$. Then, as $\ep\to0$,
$$
{\bf P}\{{}\|X\|\le \ep\}\sim {\cal C}(X)\cdot
\ep^{\gamma} \exp\left(-\ \frac {2\ell-1}{2}
\left(\frac {\th_{\ell}}{2\ell\sin\frac {\pi}{2\ell}}\right)
^{\frac {2\ell}{2\ell-1}}\ep^{-\frac {2}{2\ell-1}} \right).
\eqno(1.7)$$

\noindent Here we denote

$$\gamma=-\ell+\frac {\varkappa+1}{2\ell-1},\qquad
\th_{\ell}=\int\limits_0^1p_{\ell}^{-1/(2\ell)}(x)\,dx,$$

\noindent and the constant ${\cal C}(X)$ is given by
$${\cal C}(X)=C_{\rm dist}(X)\cdot\frac {(2\pi)^{\ell/2}
\left(\pi/\th_{\ell}\right)
^{\ell\gamma}\left(\sin \frac {\pi}{2\ell}\right)
^{\frac {1 +\gamma}{2}}} {(2\ell-1)^{1/2}
\left(\frac {\pi}{2\ell}\right)^{1+\frac {\gamma}2}
{\Gamma\vphantom)}^{\ell}\left(\ell-\frac {\varkappa}{2\ell}\right)},$$
where $C_{\rm dist}(X)$ is the so-called {\bf distortion constant}
$$C_{\rm dist}(X)\equiv\prod_{n=1}^{\infty}
\frac {\mu_n^{1/2}}{\left(\pi/\th_{\ell} \cdot
\left[n+\ell-1-\frac {\varkappa}{2\ell}\right]\right)^{\ell}},
\eqno(1.8)$$

\noindent and $\mu_n=(\lambda_n(X))^{-1}$ are the eigenvalues of the problem
(1.3)}.
\medskip

\noindent {\bf Proof}. According to comparison principle \cite{Li} and
to formulas (1.4), taking into account Remark 1 we have
$$ {\bf P}\{{}\|X\|\le \ep \}=
{\bf P}\left\{\sum_{n=1}^\infty \lambda_n \eta_n^2\le\ep ^2 \right\}=
{\bf P}\left\{\sum_{n=1}^\infty \frac {{\eta'}_n^2}{\mu'_n}+
\sum_{n=1}^\infty \frac {{\eta''}_n^2}{\mu''_n}\le\ep ^2 \right\}=$$
$$=\prod_{n=1}^{\infty}\left(\frac {\mu'_n}{\widetilde\mu'_n}\cdot
\frac {\mu''_n}{\widetilde\mu''_n}\right)^{1/2} \cdot
{\bf P}\left\{V'+V''\le\ep ^2 \right\},
\eqno(1.9)$$
where
$$V'=\sum_{n=1}^\infty \frac {{\eta'}_n^2}{\widetilde\mu'_n},\qquad
V''=\sum_{n=1}^\infty \frac {{\eta''}_n^2}{\widetilde\mu''_n},\eqno(1.10)$$
$\eta'_n$ and $\eta''_n$ are two independent sequences of independent
standard Gaussian r.v.'s,
$$\widetilde\mu'_n=\Bigl(\frac {2\pi n+\rho'}{\th_{\ell}}\Bigr)^{2\ell},
\qquad \widetilde\mu''_n=
\Bigl(\frac {2\pi n+\rho''}{\th_{\ell}}\Bigr)^{2\ell}.\eqno(1.11)$$

The asymptotic behavior of small ball probabilities for the infinite sums
(1.10) with coefficients of the form (1.11) was derived in
\cite[Theorem 6.2]{NN}. The asymptotics of
${\bf P}\left\{V'+V''\le\ep ^2 \right\}$ as $\ep\to0$ can be deduced from the
asymptotics of ${\bf P}\left\{V'\le\ep ^2 \right\}$ and
${\bf P}\left\{V''\le\ep ^2 \right\}$ by Lemma 0.1. Finally, the infinite
product in (1.9) differs from (1.8) by the multiplier which converges due to
(1.5). After simplification, we arrive at (1.7).\hfill$\square$\medskip

\section{Slepian process and related processes}

Consider the generalized Slepian process $S^{(c)}$, that is a stationary zero
mean Gaussian process with covariance function
$$G_{S^{(c)}}(t,s)= c- |t-s|, \qquad t,s\in [0,1].
$$
It is easy to check that $G_{S^{(c)}}$ is indeed a covariance for $c\ge1/2$.
Remark that for $c\ge1$ we have the distributional equality
$$S^{(c)}(t)\stackrel {d}{=}W(t+c)-W(t),\qquad 0\le t\le1,$$
where $W(t)$ is a standard Wiener process. The conventional Slepian process
\cite{Sl} corresponds to $c=1$; the small ball asymptotics for this process
was derived in \cite{NO}.

The direct calculation shows, see \cite{GL}, that $G_{S^{(c)}}$ is the Green
function of the BVP
$${\cal L}_{S^{(c)}}u\equiv -{\textstyle\frac 12} u''=\mu u\ \
\mbox{\rm on} \ \ [0,1],\qquad
u'(0)+u'(1)=0,\ \ (2c-1)u'(0)-(u(0)+u(1))=0.\eqno(2.1)$$

{\bf Theorem 2.1}. \ {\bf 1}. {\it Let $c=1/2$. Then, as $\ep\to0$,}
$${\bf P}\{\|S^{(1/2)}\|\le\ep\}\sim\frac {4\,\ep}{\sqrt{\pi}}\cdot
\exp(-{\textstyle\frac 14}\ep^{-2}).\eqno(2.2)$$
{\bf 2}. {\it Let $c>1/2$. Then, as $\ep\to0$,}
$${\bf P}\{\|S^{(c)}\|\le\ep\}\sim\frac {4\sqrt{2}\,\ep^2}{\sqrt{\pi(2c-1)}}
\cdot\exp(-{\textstyle\frac 14}\ep^{-2}).\eqno(2.3)$$

{\bf Remark 3}. When $c>1/2$, the problem (2.1) does not satisfy assumptions
of Theorem 1.2, since $\varkappa=2=2\ell^2$. This is only a formal
difficulty, and the ways to avoid it are well known, see, e.g.,
\cite[Proposition 6.4]{NN} and \cite{NO}. In this particular case, however,
it is more simple to use the available results.\medskip

\noindent {\bf Proof}. {\bf 1}. Set $\zeta=\sqrt{2\mu}$. Substituting the
general solution of the equation (2.1)
$u(t)=c_1\sin(\zeta t)+c_2\cos(\zeta t)$ into boundary conditions we deduce
that $\mu_n=\frac 12r_n^2$, where $r_1<r_2<\ldots$ are positive roots of the
equation
$$F^{(c)}(\zeta)\equiv 2+2\cos(\zeta)-(2c-1)\zeta\sin(\zeta)=0.$$
It is easy to see that for $c=1/2$ the spectrum of the problem (2.1) consists
of double eigenvalues $\mu'_n=\mu''_n=2(\pi n-\frac {\pi}2)^2$,
$n\in{\mathbb N}$. Thus, we have the distributional equality
$$\|S^{(1/2)}\|^2\stackrel {d}{=}\frac 12(\|W_1\|^2+\|W_2\|^2),$$
where $W_1$ and $W_2$ are independent standard Wiener processes.
The asymptotics of ${\bf P}\{\|W\|\le\ep\}$ as $\ep\to0$ is well known.
Applying Lemma 0.1 we arrive at (2.2).\medskip

{\bf 2}. When $c>1/2$ it is obvious that
$\frac {|F^{(c)}(0)|}{|F^{(1)}(0)|}=1$ and
$\frac {|F^{(c)}(\zeta)|}{|F^{(1)}(\zeta)|}\rightrightarrows 2c-1$ as
$|\zeta|=\pi(N+\frac 12)$, $N\to\infty$. Let us apply the comparison theorem
\cite{Li} to the processes $S^{(c)}$ and $S^{(1)}$. Then we apply Jensen's
Theorem, see \cite[\S3.6.1]{T}, to the functions $F^{(c)}$ and $F^{(1)}$ and
obtain
$${\bf P}\{\|S^{(c)}\|\le\ep\}\sim\frac 1{\sqrt{2c-1}}\cdot
{\bf P}\{\|S^{(1)}\|\le\ep\}, \qquad  \ep\to0.$$
The asymptotics of the last expression, as we mentioned, was established in
\cite{NO}. This gives (2.3).\hfill$\square$\medskip

Now we consider $m$-times integrated process
$(S^{(c)})_m^{[\beta_1,\,...,\,\beta_m]}(t)$. Following \cite{Na} we
introduce the notation
$$\widetilde\ep_{\ell}=\left(\ep\sqrt {\ell\sin\frac {\pi}{2\ell}}\right)
^{\frac {1}{2\ell-1}};\qquad {\mathfrak D}_{\ell}=\frac
{2\ell-1}{2\ell\sin\frac {\pi}{2\ell}};\eqno(2.4)$$
for $j=1,\dots,m$
$$k_j=\ \left\{
\begin{array}
{lll} m-j,& \mbox {if}& \beta_j=0,\\
m+1+j,& \mbox {if}&\beta_j=1,
\end{array}\right.
\qquad k'_j=2m+1-k_j.
\eqno(2.5)$$

{\bf Theorem 2.2}. {\it Let $m\in{\mathbb N}$. Then, as $\ep\to0$,}\medskip

{\bf 1}. {\it for $c=1/2$}
\begin{multline*}
{\bf P}\{\|(S^{(1/2)})_m^{[\beta_1,\,...,\,\beta_m]}\|\le\ep\}\sim\\
\sim\frac {(2m+2)^{\frac m2+1}}
{|{\mathfrak V}(z_{m+1}^{k_1},z_{m+1}^{k_2},
\dots,z_{m+1}^{k_m})|\cdot \sqrt{\prod \limits_{j=1}^m |1+z_{m+1}^{k_j}|^2+
\prod \limits_{j=1}^m |1+z_{m+1}^{k'_j}|^2}}\cdot\\
\cdot\frac {2\widetilde\ep_{m+1}}
{\sqrt {\pi {\mathfrak D}_{m+1}}}
\exp\left(-\ \frac {{\mathfrak D}_{m+1}}{2\widetilde\ep_{m+1}^2}
\right).\tag{2.6}
\end{multline*}

{\bf 2}. {\it for $c>1/2$}
\begin{multline*}
{\bf P}\{\|(S^{(c)})_m^{[\beta_1,\,...,\,\beta_m]}\|\le\ep\}\sim
\frac {(2m+2)^{\frac {m+1}2}\sqrt {2\sin\frac {\pi\vphantom2}{2m+2}}}
{|{\mathfrak V}(z_{m+1}^{k_1},z_{m+1}^{k_2},\dots,z_{m+1}^{k_m})|
\sqrt {2c-1}}\cdot \\
\cdot\frac {2\widetilde\ep_{m+1}^2}
{\sqrt {\pi {\mathfrak D}_{m+1}}}
\exp\left(-\ \frac {{\mathfrak D}_{m+1}}{2\widetilde\ep_{m+1}^2}
\right).\tag{2.7}
\end{multline*}

\noindent {\bf Proof}. The BVP generated by integrated process can be
expressed in terms of BVP generated by original process due to
\cite[Theorem 2.1]{NN}. Applying this theorem to (2.1) we get
$$
\left\{
\begin{aligned}
&{\cal L}_{S^{(c)}_m}u\equiv (-1)^{m+1}\cdot
{\textstyle\frac 12} u^{(2m+2)}=\mu u \quad \text{on} \quad [0,1],\\
&u(\beta_m)=u'(\beta_{m-1})=\dots=u^{(m-1)}(\beta_1)=0,\\
&u^{(m+1)}(0)+u^{(m+1)}(1)=0,\quad
(2c-1)u^{(m+1)}(0)-(u^{(m)}(0)+u^{(m)}(1))=0,\\
&u^{(m+2)}(1-\beta_1)=u^{(m+3)}(1-\beta_2)=\dots=u^{(2m+1)}(1-\beta_m)=0.\\
\end{aligned}
\right.
\eqno(2.8)$$
Using (2.5) we can rewrite the second and the fourth rows in (2.8) as
follows:
$$u^{(k_j)}(0)=u^{(k'_j)}(1)=0,\qquad j=1,\dots,m.$$

Since the BVP (2.8) satisfies all assumptions of Theorem 1.2 (with
$\ell=m+1$), to prove (2.6)-(2.7) we only need to calculate the distortion
constants (1.8). Note that $\th_{\ell}=2^{\frac 1{2m+2}}$.\medskip

Put $\zeta=(2\mu)^{\frac 1{2m+2}}$. Then the general solution of the equation
(2.8) is
$$u(t)=\sum\limits_{j=0}^{2m+1}c_j\exp(\omega_j\zeta t),\eqno(2.9)$$
where $\omega_j=z_{m+1}^j$.\medskip

{\bf 1}. If $c=\frac 12$ then $\varkappa=(2m+1)(m+1)$. Substituting (2.9)
into boundary conditions we deduce that $\mu_n=\frac 12r_n^{2m+2}$, where
$r_1<r_2<\ldots$ are positive roots of the entire function
\begin{multline*}
{\cal F}(\zeta)\equiv\\
\det
\begin{bmatrix}
1&\omega_1^{k_1}&\dots&\omega_m^{k_1}&(-1)^{k_1}&\dots&(-\omega_m)^{k_1}\\
\hdotsfor[3]{7} \\
1&\omega_1^{k_m}&\dots&\omega_m^{k_m}&(-1)^{k_m}&\dots&(-\omega_m)^{k_m}\\
1+e^{i\zeta}&\omega_1^m(1+e^{i\omega_1\zeta})&\dots&
\dots&(-1)^m(1+e^{-i\zeta})&\dots&\dots\\
1+e^{i\zeta}&\omega_1^{m+1}(1+e^{i\omega_1\zeta})&\dots&
\dots&(-1)^{m+1}(1+e^{-i\zeta})&\dots&\dots\\
e^{i\zeta}&\omega_1^{k'_1}e^{i\omega_1\zeta}&\dots&
\omega_m^{k'_1}e^{i\omega_m\zeta}&(-1)^{k'_1}e^{-i\zeta}&\dots&
(-\omega_m)^{k'_1}e^{-i\omega_m\zeta}\\
\hdotsfor[3]{7} \\
e^{i\zeta}&\omega_1^{k'_m}e^{i\omega_1\zeta}&\dots&
\omega_m^{k'_m}e^{i\omega_m\zeta}&(-1)^{k'_m}e^{-i\zeta}&\dots&
(-\omega_m)^{k'_m}e^{-i\omega_m\zeta}\\
\end{bmatrix}.
\end{multline*}
Therefore,
$$C_{\rm dist}(S^{(1/2)}_m)=\prod_{n=1}^{\infty}\biggl(\frac {r_n}
{\pi n-\frac {\pi}2}\biggr)^{m+1}.
$$

Since $|{\cal F}(\zeta)|\equiv |{\cal F}(\omega_1\zeta)|$, the set of all
nonzero roots of the function ${\cal F}$ consists of $2m+2$ sequences
$\omega_jr_n$, $j=0,\dots,2m+1$, $n\in\mathbb N$.

According to \cite[\S4, Theorem 2]{Nm}, the relation
$${\cal F}(\zeta)=\exp(-i\omega_1\zeta)\exp(-i\omega_2\zeta)\dots
\exp(-i\omega_m\zeta)\cdot\Bigl(\Phi(\zeta)+O(|\zeta|^{-1})\Bigr),
\eqno(2.10)$$
holds true for $|\zeta|\to\infty$ and $|\arg(\zeta)|\le \frac {\pi}{2m+2}$.
Here
$$\Phi(\zeta)=\det
\begin{bmatrix}
1&\omega_1^{k_1}&\dots&\omega_m^{k_1}&(-1)^{k_1}&0&\dots&0\\
\hdotsfor[3]{8} \\
1&\omega_1^{k_m}&\dots&\omega_m^{k_m}&(-1)^{k_m}&0&\dots&0\\
1+e^{i\zeta}&\omega_1^m&\dots&\omega_m^m&(-1)^m(1+e^{-i\zeta})
&(-\omega_1)^m&\dots&(-\omega_m)^m\\
1+e^{i\zeta}&\omega_1^{m+1}&\dots&\omega_m^{m+1}&(-1)^{m+1}(1+e^{-i\zeta})
&(-\omega_1)^{m+1}&\dots&(-\omega_m)^{m+1}\\
e^{i\zeta}&0&\dots&0&(-1)^{k'_1}e^{-i\zeta}&(-\omega_1)^{k'_1}&\dots&
(-\omega_m)^{k'_1}\\
\hdotsfor[3]{8} \\
e^{i\zeta}&0&\dots&0&(-1)^{k'_m}e^{-i\zeta}&(-\omega_1)^{k'_m}&\dots&
(-\omega_m)^{k'_m}\\
\end{bmatrix}.$$

Expanding this determinant in the elements of the first and the $(m+2)$-nd
columns we obtain
$$|\Phi(\zeta)|={\cal M}\cdot |\exp(i\zeta)+\exp(-i\zeta)+R|,
\eqno(2.11)$$
where
\begin{multline*}
{\cal M}=
|{\mathfrak V}(\omega_1^{k_1},\dots,\omega_1^{k_m},\omega_1^m)\cdot
{\mathfrak V}(\omega_1^{m+1},\omega_1^{k'_1},\dots,\omega_1^{k'_m})+\\
+{\mathfrak V}(\omega_1^{k_1},\dots,\omega_1^{k_m},\omega_1^{m+1})\cdot
{\mathfrak V}(\omega_1^m,\omega_1^{k'_1},\dots,\omega_1^{k'_m})|,
\end{multline*}
while $R$ is a constant which is inessential for us.

Following \cite{Na}, for arbitrary $\delta>-1$ we introduce the function
$$\Psi_{\delta}(\zeta)=\psi_{\delta}(\zeta)
\psi_{\delta}(\omega_1\zeta)\psi_{\delta}(\omega_2\zeta)\dots
\psi_{\delta}(\omega_m\zeta),\eqno(2.12)$$
where
$$\psi_{\delta}(\zeta)=\frac {\Gamma^2(1+\delta)}
{\Gamma\left(1+\delta+\frac {\zeta}{\pi}\right)
\Gamma\left(1+\delta-\frac {\zeta}{\pi}\right)}
=\prod\limits_{n=1}^{\infty}
\left(1-\frac {\zeta^2}{(\pi(n+\delta))^2}\right).$$
It is easy to check, see \cite[Lemma 1.3]{Na}), that
$$\psi_{\delta}(\zeta)\sim \Gamma^2(1+\delta)\pi^{2\delta}
\zeta^{-2\delta-1}\cos(\zeta-\pi(\delta+1/2)),\eqno(2.13)$$
as $\zeta\to\infty$, $|\arg(\zeta)|\le \phi_0<\pi$. Moreover, the convergence
is uniform in $\arg(\zeta)$.

Setting $\delta=-1/2$ we deduce from (2.10)-(2.13)
$$\frac {|{\cal F}(\zeta)|}{|\Psi_{\delta}(\zeta)|}\to
2^{m+1}{\cal M},\eqno(2.14)$$
as $\zeta\to\infty$, $\arg(\zeta)\ne \frac{\pi j}{2m+2}$, $j\in\mathbb Z$.

By Jensen's Theorem
$$C^2_{\rm dist}(S^{(1/2)}_m)=
\frac {|{\cal F}(0)|}{|\Psi_{\delta}(0)|}\cdot \exp\biggl\{
\lim_{\rho \to \infty}\frac 1{2\pi} \int_0^{2\pi}
\ln \frac {|\Psi_{\delta}(\rho e^{i\theta})|}
{|{\cal F}(\rho e^{i\theta})|}\,d\theta \biggr\}.$$
The integrand obviously has a summable majorant. The Lebesgue Dominated
Convergence Theorem gives, in view of (2.14),
$$C^2_{\rm dist}(S^{(1/2)}_m)=
\frac {|{\cal F}(0)|}{2^{m+1}{\cal M}}=
\frac {4|{\mathfrak V}(1,\omega_1,\omega_1^2,\dots,\omega_1^{2m+1})|}
{2^{m+1}{\cal M}}=\frac {4(m+1)^{m+1}}{\cal M}.$$
It remains to take into account that, due to (2.5),
$${\cal M}=
|{\mathfrak V}(\omega_1^{k_1},\dots,\omega_1^{k_m})|^2\cdot
\Big(\prod \limits_{j=1}^m |1+\omega_1^{k_j}|^2+
\prod \limits_{j=1}^m |1+\omega_1^{k'_j}|^2\Big).$$
After some simplification we arrive at (2.6).\medskip

{\bf 2}. If $c>\frac 12$ then $\varkappa=(2m+1)(m+1)+1$. Substituting (2.9)
into boundary conditions we deduce that
$$C_{\rm dist}(S^{(c)}_m)=\prod_{n=1}^{\infty}\biggl(\frac {r_n}
{\pi (n-\frac {m+2}{2m+2})}\biggr)^{m+1},
$$
where $r_1<r_2<\ldots$ are positive roots of the function
${\cal F}_{(c)}(\zeta)$, which arises if we change in the determinant
${\cal F}(\zeta)$ the row
$$
\begin{bmatrix}
1+e^{i\zeta}&\omega_1^m(1+e^{i\omega_1\zeta})&\dots&
\omega_m^m(1+e^{i\omega_m\zeta})
&(-1)^m(1+e^{-i\zeta})&\dots&(-\omega_m)^m(1+e^{-i\omega_m\zeta})\\
\end{bmatrix}
$$
by the row
$$
\begin{bmatrix}
1-\frac {\tau}{i\zeta}(1+e^{i\zeta})&
\omega_1^{m+1}(1-\frac {\tau}{i\omega_1\zeta}(1+e^{i\omega_1\zeta}))&\dots&
\dots&(-1)^{m+1}(1-\frac {\tau}{-i\zeta}(1+e^{-i\zeta}))&\dots&\dots\\
\end{bmatrix}
$$
(here $\tau=\frac 1{2c-1}$).

Similarly to part {\bf 1}, the relation
$${\cal F}_{(c)}(\zeta)=\exp(-i\omega_1\zeta)\exp(-i\omega_2\zeta)\dots
\exp(-i\omega_m\zeta)\cdot\Bigl(\Phi(\zeta)+O(|\zeta|^{-1})\Bigr),
$$
holds true as $|\zeta|\to\infty$ and $|\arg(\zeta)|\le \frac {\pi}{2m+2}$.
Here
$$\Phi(\zeta)=\det
\begin{bmatrix}
1&\omega_1^{k_1}&\dots&\omega_m^{k_1}&(-1)^{k_1}&0&\dots&0\\
\hdotsfor[3]{8} \\
1&\omega_1^{k_m}&\dots&\omega_m^{k_m}&(-1)^{k_m}&0&\dots&0\\
1&\omega_1^{m+1}&\dots&\omega_m^{m+1}&(-1)^{m+1}&0&\dots&0\\
1+e^{i\zeta}&\omega_1^{m+1}&\dots&\omega_m^{m+1}&(-1)^{m+1}(1+e^{-i\zeta})
&(-\omega_1)^{m+1}&\dots&(-\omega_m)^{m+1}\\
e^{i\zeta}&0&\dots&0&(-1)^{k'_1}e^{-i\zeta}&(-\omega_1)^{k'_1}&\dots&
(-\omega_m)^{k'_1}\\
\hdotsfor[3]{8} \\
e^{i\zeta}&0&\dots&0&(-1)^{k'_m}e^{-i\zeta}&(-\omega_1)^{k'_m}&\dots&
(-\omega_m)^{k'_m}\\
\end{bmatrix}.$$
Let us subtract $(m+1)$-st row from $(m+2)$-nd one. Then, expanding the
determinant in the elements of the first column we obtain
\begin{multline*}
|\Phi(\zeta)|=
|{\mathfrak V}(\omega_1^{k_1},\dots,\omega_1^{k_m},\omega_1^{m+1})| \cdot
|{\mathfrak V}(\omega_1^{m+1},\omega_1^{k'_1},\dots,\omega_1^{k'_m})| \cdot\\
\cdot|\exp(-i\zeta)-\omega_1^\varkappa\exp(i\zeta)|=
|{\mathfrak V}(\omega_1^{k_1},\dots,\omega_1^{k_m})|\cdot
|{\mathfrak V}(\omega_1^{k'_1},\dots,\omega_1^{k'_m})|\cdot\\
\cdot\prod\limits_{j=1}^{m}
|(\omega_1^{m+1}-\omega_1^{k_j})(\omega_1^{m+1}-\omega_1^{k'_j})|\cdot
|\exp(-i\zeta)+\omega_1\exp(i\zeta)|=\\
=|{\mathfrak V}(\omega_1^{k_1},\dots,\omega_1^{k_m})|^2\cdot
\frac {2m+2}{|1-\omega_1|}\cdot
|\exp(-i\zeta)+\omega_1\exp(i\zeta)|.
\end{multline*}

Hence, setting $\delta=-\frac {m+2}{2m+2}$ and taking into account
$|{\cal F}_{(c)}(\zeta)|\equiv |{\cal F}_{(c)}(\omega_1\zeta)|$, we obtain
$$\frac {|\zeta {\cal F}_{(c)}(\zeta)|}{|\Psi_{\delta}(\zeta)|}\to
\frac {2^{m+1}|{\mathfrak V}(\omega_1^{k_1},\dots,\omega_1^{k_m})|^2 }
{\Gamma^{2m+2}(1+\delta)\pi^{(2m+2)\delta}}\cdot\frac {2m+2}{|1-\omega_1|},
$$
as $\zeta\to\infty$, $\arg(\zeta)\ne \frac{\pi j}{2m+2}$, $j\in\mathbb Z$.

Applying Jensen's Theorem similarly to part {\bf 1}, we arrive at
\begin{multline*}
C^2_{\rm dist}(S^{(c)}_m)=
\frac {\Gamma^{2m+2}(1+\delta)\pi^{(2m+2)\delta}|1-\omega_1|}
{2^{m+1}|{\mathfrak V}(\omega_1^{k_1},\dots,\omega_1^{k_m})|^2(2m+2)}
\cdot\big|\zeta {\cal F}_{(c)}(\zeta)\big|_{\zeta=0}=\\
=\frac {2\tau(m+1)^m\Gamma^{2m+2}(1+\delta)|1-\omega_1|}
{\pi^{m+2}|{\mathfrak V}(\omega_1^{k_1},\dots,\omega_1^{k_m})|^2}.
\end{multline*}
Since $|1-\omega_1|=2\sin\frac {\pi}{2m+2}$, this gives (2.7) after some
simplification.\hfill$\square$\medskip

{\bf Remark 4}. Let $c>1/2$. Using (2.7) and the extremal properties of the
Vandermond determinants \cite{N2} one can see that among all $m$-times
integrated processes $(S^{(c)})_m^{[\beta_1,\,...,\,\beta_m]}$ the processes
$(S^{(c)})_m^{[0,\,...,\,0]}$ and $(S^{(c)})_m^{[1,\,...,\,1]}$ have the
largest small ball constant while the {\bf Euler integrated} processes
$(S^{(c)})_m^{[0,1,0,...]}$ and $(S^{(c)})_m^{[1,0,1,...]}$ have the smallest
one. We conjecture that this is true also for $c=1/2$ but this problem is
open yet.\medskip

We also point out a curious relation arising when one compares the small ball
asymptotics for the process $(S^{(c)})_m$
, $c>1/2$, and for integrated Ornstein -- Uhlenbeck process, see
\cite[Theorem 2.2]{Na}. Since corresponding BVPs have the same parameters
$\varkappa$, these asymptotics differ only by a constant. Rather unexpected
is the fact that this constant equals $\sqrt{\frac{e}{2(2c-1)}}$ and
therefore depends neither on $\beta_j$ nor even on $m$.

\section{Integrated centered Brownian bridge and \\ related processes}

The most famous process generating the BVP with non-separated boundary
conditions is the centered Brownian bridge ${\overline B}(t)$; its spectrum
was derived for the first time in \cite{W}. Note that this BVP
$${\cal L}_{\overline B}u\equiv - u''=\mu u\quad \text{on} \quad [0,1],\qquad
u(0)-u(1)=0,\quad u'(0)-u'(1)=0\eqno(3.1)$$
has a zero eigenvalue with constant eigenfunction $\varphi_0(t)\equiv1$.
Hence the covariance $G_{\overline B}(t,s)$ is the generalized Green function
of the problem (3.1). We show later that this is a typical situation for
centered processes. In this case Theorem 2.1 \cite{NN} is not applicable. To
study integrated processes we need two auxiliary statements.\medskip

{\bf Theorem 3.1}. \ {\bf 1}. { \it Let the BVP (1.3) have a zero eigenvalue
with constant eigenfunction $\varphi_0(t)\equiv1$. Let the kernel $G(t,s)$
be the generalized Green function of the problem (1.3). Then the integrated
kernel
$${\cal G}_1(t,s)=\int_0^t\int_0^s G(x,y)\,dxdy\eqno(3.2)$$

\noindent is (conventional) Green function of the BVP
$${\cal L}_1u\equiv-({\cal L}u')'=\mu u\quad \mbox
{\rm on}\quad [0,1], \qquad u\in {\cal D}({\cal L}_1),\eqno(3.3)$$

\noindent where the domain ${\cal D}({\cal L}_1)$ consists of the functions
$u\in W_2^{2\ell+2}(0,1)$ satisfying the boundary conditions
$$u(0)=0;\qquad u(1)=0;\qquad u'\in {\cal D}({\cal L}).\eqno(3.4)$$

{\bf 2}. Let the kernel $G(t,s)$ be the Green function of the problem
(3.3)-(3.4). Then the centered kernel
$${\overline G}(t,s)=G(t,s)-g(t)-g(s)+\overline g \eqno(3.5)$$
(here $g(t)=\int\limits_0^1G(t,s)\,ds$,
$\overline g=\int\limits_0^1g(t)\,dt$) is the generalized Green function of
the BVP \footnote {The problem (3.6)-(3.7) obviously has a zero eigenvalue
with constant eigenfunction $\varphi_0(t)\equiv1$.}
$$\widetilde {\cal L}_1u\equiv-({\cal L}u')'=\mu u\quad \mbox
{\rm on}\quad [0,1], \qquad u\in {\cal D}(\widetilde {\cal L}_1),\eqno(3.6)$$
\noindent where the domain ${\cal D}(\widetilde {\cal L}_1)$ consists of the
functions $u\in W_2^{2\ell+2}(0,1)$ satisfying the boundary conditions}
$$u(0)-u(1)=0;\qquad u'\in {\cal D}({\cal L}); \qquad
({\cal L}u')(0)-({\cal L}u')(1)=0.\eqno(3.7)$$

{\bf Remark 5}. Let $X(t)$, $0\le t\le 1$ be a Gaussian process with zero
mean. It is well known that the covariance of the integrated process
$G_{X_1^{[0]}}$ can be expressed in terms of the original covariance $G_X$ by
formula (3.2). Note that under assumptions of part {\bf 1} of our Theorem the
processes $X_1^{[0]}$ and $X_1^{[1]}$ coincide almost surely. It is easy to
show that the covariance of the centered process $G_{\overline X}$ can be
expressed in terms of the original covarisnce $G_X$ by formula (3.5).\medskip

{\bf Remark 6}. It is easy to see that the differential expression (1.1) can
be represented in the form (3.3) iff $p_0\equiv0$.\medskip

\noindent {\bf Proof}. {\bf 1}. The first boundary condition in (3.4) is
trivially satisfied for the function ${\cal G}_1$ while the second one is
satisfied due to (0.1). Further, differentiating (3.2) w.r.t. $t$ we obtain
$$({\cal G}_1)'_t(t,s)=\int_0^s G(t,y)\,dy,$$
whence the other boundary conditions follow by linearity of the set
${\cal D}({\cal L})$. Since ${\cal L}G(t,s)=\delta(t-s)-1$, we obtain
consequently
$${\cal L}({\cal G}_1)'_t(t,s) =\chi_{\mathbb R_+}(t-s)-s;\qquad
({\cal L}({\cal G}_1)'_t)'_t=-\delta(t-s),$$
and the statement follows.\medskip

{\bf 2}. The orthogonality condition (0.1) follows from the definition of
${\overline G}$:
$$\int_0^1{\overline G}(t,s)\,ds=g(s)-g(s)-\overline g+\overline g=0.$$
The first two conditions in (3.4) provide
$${\overline G}(0,s)=-g(s)+\overline g={\overline G}(1,s).$$
Differentiating (3.5) w.r.t. $t$ we obtain
$${\overline G}{\vphantom G}'_t(t,s)=G'_t(t,s)-\int_0^1 G'_t(t,y)\,dy,$$
that gives ${\overline G}{\vphantom G}'_t\in{\cal D}({\cal L})$.
Since ${\cal L}_1G(t,s)=\delta(t-s)$, we have
$$\widetilde {\cal L}_1{\overline G}(t,s)=\delta(t-s)-
\int_0^1 \delta(t-y)\,dy=\delta(t-s)-1.$$
Finally, the last boundary condition follows from
$$({\cal L}{\overline G}{\vphantom G}'_t)(0,s)-
({\cal L}{\overline G}{\vphantom G}'_t)(1,s)=
\int_0^1\widetilde {\cal L}_1{\overline G}(t,s)\,dt=
\int_0^1(\delta(t-s)-1)\,dt=0,$$
and the second statement is also proved.\hfill$\square$\medskip

Now we define the sequence of integrated centered analogues of Brownian
bridge. We set
$$B_{\{0\}}(t)=B(t); \qquad
B_{\{l\}}(t)=\int_0^t\overline {B_{\{l-1\}}}(s)\,ds, \quad
l\in\mathbb N.$$
Theorem 3.1 allows us to write down the BVPs generated by processes
$B_{\{l\}}$ and $\overline {B_{\{l\}}}$. We begin from the second process
because its eigenvalues can be derived explicitly. This permits us to derive
the small ball asymptotics without using Theorem 1.2.\medskip

{\bf Theorem 3.2}. { \it Let $l\in{\mathbb N}_0$. Then, as $\ep\to0$,
$${\bf P}\{\|\overline {B_{\{l\}}}\|\le\ep\}\sim
\sqrt {2l+2}\cdot\frac {\ep_{l+1}^{-(2l+1)}}
{\sqrt {\pi {\mathfrak D}_{l+1}}}
\exp\left(-\ \frac {{\mathfrak D}_{l+1}}{2\ep_{l+1}^2}\right),
\eqno(3.8)$$
where $\ep_{\ell}=\left(\ep\sqrt {2\ell\sin\frac {\pi}{2\ell}}\right)
^{\frac {1}{2\ell-1}}$ while the quantity ${\mathfrak D}_{\ell}$ was defined
in (2.4).}\medskip

{\bf Remark 7}. The multiplier before the exponent in (3.8) equals
$\sqrt {\frac {2l+2}{2l+1}}\cdot\frac {\ep^{-1}}{\sqrt {\pi}}$. For $l=0$
(3.8) coincides with formula obtained in \cite[\S3]{BNO}.\medskip

\noindent {\bf Proof}. Applying $l$ times in turns the first and the second
statements of Theorem 3.1 to the problem (3.1) we deduce that the covariance
$G_{\overline {B_{\{l\}}}}(t,s)$ is the generalized Green function of the BVP
with periodic boundary conditions
$${\cal L}_{\overline {B_{\{l\}}}}u\equiv  (-1)^{l+1}u^{(2l+2)}=\mu u \quad
\mbox{\rm on} \quad [0,1],\qquad
u^{(j)}(0)-u^{(j)}(1)=0,\quad j=0,1,\dots,2l+1.$$
Whence the operator ${\cal L}_{\overline {B_{\{l\}}}}$ coincides with
$({\cal L}_{\overline B})^{l+1}$. Therefore, its spectrum is double,
excluding zero eigenvalue which is inessential for us due to the
orthogonality condition (0.1): $\mu'_n=\mu''_n=(2\pi n)^{2l+2}$. Thus, we
have distributional equality
$$\|\overline {B_{\{l\}}}\|^2\stackrel {d}{=}
\sum_{n=1}^\infty \frac {{\eta'}_n^2}{(2\pi n)^{2l+2}}+
\sum_{n=1}^\infty \frac {{\eta''}_n^2}{(2\pi n)^{2l+2}},$$
where $\eta'_n$ and $\eta''_n$ are two independent sequences of independent
standard Gaussian r.v.'s. Using \cite[Theorem 6.2]{NN} and Lemma 0.1 we
arrive at (3.8).\hfill$\square$\medskip

{\bf Theorem 3.3}. {\it Let $l\in{\mathbb N}_0$. Then, as $\ep\to0$,
$${\bf P}\{\|B_{\{l\}}\|\le\ep\}\sim
(2l+2)\sqrt {\textstyle\sin\frac {\pi\vphantom2}{2l+2}}\cdot
\frac {\ep_{l+1}^{-2l}}{\sqrt {\pi {\mathfrak D}_{l+1}}}
\exp\left(-\ \frac {{\mathfrak D}_{l+1}}{2\ep_{l+1}^2}\right),
\eqno(3.9)$$
with the same notation as in Theorem 3.2}.\medskip

\noindent {\bf Proof}. Similarly to Theorem 3.2, the covariance
$G_{B_{\{l\}}}(t,s)$ is the Green function of the BVP
$$\left\{
\begin{aligned}
&{\cal L}_{B_{\{l\}}}u\equiv  (-1)^{l+1}u^{(2l+2)}=\mu u\quad \text{on}
\quad [0,1],\\
&u(0)=u(1)=0,\quad u^{(j)}(0)-u^{(j)}(1)=0,\ \ j=1,\dots,2l.\\
\end{aligned}
\right.\eqno(3.10)$$
Since the problem (3.10) satisfies all assumptions of Theorem 1.2 (with
$\ell=l+1$), to prove (3.9) we only need to calculate the distortion
constant. Note that $\th_{\ell}=1$ and $\varkappa=(2l+1)l$.

Put $\zeta=\mu^{\frac 1{2l+2}}$. Substituting the general solution of the
equation (3.10) into boundary conditions, we deduce that
$$C_{\rm dist}(B_{\{l\}})=\prod_{n=1}^{\infty}\biggl(\frac {r_n}
{\pi (n+\frac l{2l+2})}\biggr)^{l+1},
$$
where $r_1<r_2<\ldots$ are positive roots of the entire function
$${\mathfrak F}(\zeta)\equiv\det
\begin{bmatrix}
1&1&\dots&1&1&\dots&1\\
e^{i\zeta}&e^{i\omega_1\zeta}&\dots&
e^{i\omega_l\zeta}&e^{-i\zeta}&\dots&e^{-i\omega_l\zeta}\\
1-e^{i\zeta}&\omega_1(1-e^{i\omega_1\zeta})&\dots&
\dots&(-1)(1-e^{-i\zeta})&\dots&\dots\\
1-e^{i\zeta}&\omega_1^2(1-e^{i\omega_1\zeta})&\dots&
\dots&(-1)^2(1-e^{-i\zeta})&\dots&\dots\\
\hdotsfor[3]{7} \\
1-e^{i\zeta}&\omega_1^{2l}(1-e^{i\omega_1\zeta})&\dots&
\dots&(-1)^{2l}(1-e^{-i\zeta})&\dots&\dots\\
\end{bmatrix},
$$
while $\omega_j=z_{l+1}^j$.

Subtracting the first row from the second one, similarly to Theorem 2.2 we
obtain for $|\zeta|\to\infty$ and $|\arg(\zeta)|\le \frac {\pi}{2l+2}$
$${\mathfrak F}(\zeta)=(-1)^{l+1}\exp(-i\omega_1\zeta)\exp(-i\omega_2\zeta)
\dots\exp(-i\omega_l\zeta)\cdot\Bigl(\Phi(\zeta)+O(|\zeta|^{-1})\Bigr),
$$

\noindent where
$$\Phi(\zeta)=\det
\begin{bmatrix}
1&1&\dots&1&1&0&\dots&0\\
1-e^{i\zeta}&1&\dots&1&1-e^{-i\zeta}
&1&\dots&1\\
1-e^{i\zeta}&\omega_1&\dots&\omega_l&(-1)(1-e^{-i\zeta})
&-\omega_1&\dots&-\omega_l\\
\hdotsfor[3]{8} \\
1-e^{i\zeta}&\omega_1^{2l}&\dots&\omega_l^{2l}&(-1)^{2l}(1-e^{-i\zeta})
&(-\omega_1)^{2l}&\dots&(-\omega_l)^{2l}\\
\end{bmatrix}.$$

Expanding this determinant in the elements of the first row we derive
\begin{multline*}
|\Phi(\zeta)|=\frac {2|{\mathfrak V}(1,\omega_1,\dots,\omega_{2l})|}
{|1-\omega_1|}\cdot |(1-\exp(i\zeta))(\exp(-i\zeta)-\omega_1)|=\\
=\frac {2(2l+2)^l}{|1-\omega_1|}\cdot
|\exp(-i\zeta)+\omega_1\exp(i\zeta)-(1+\omega_1)|.
\end{multline*}

Setting $\delta=\frac l{2l+2}$ we obtain in view of
$|{\mathfrak F}(\zeta)|\equiv |{\mathfrak F}(\omega_1\zeta)|$
$$\frac {|{\mathfrak F}(\zeta)|}
{|\zeta^{2l+1}\prod \limits_{j=0}^{l}\psi_{\delta}(\omega_j\zeta)|}\to
\frac {2^{l+2}(2l+2)^l}{\Gamma^{2l+2}(1+\delta)\pi^l\,|1-\omega_1|},
$$
for $\zeta\to\infty$, $\arg(\zeta)\ne \frac{\pi j}{2l+2}$, $j\in\mathbb Z$.

This implies, similarly to the proof of Theorem 2.2,
$$C^2_{\rm dist}(B_{\{l\}})=
\frac {\Gamma^{2l+2}(1+\delta)\pi^l\,|1-\omega_1|}
{2^{l+2}(2l+2)^l}\cdot
\left|\frac {{\mathfrak F}(\zeta)}{\zeta^{2l+1}} \right|_{\zeta=0}.
$$
Since
$$\frac {-{\mathfrak F}(\zeta)}{\zeta^{2l+1}}=
\det
\begin{bmatrix}
1&1&\dots&1&1&\dots&1\\
\frac {1-e^{i\zeta}}{\zeta}&\frac {1-e^{i\omega_1\zeta}}{\zeta}&\dots&
\frac {1-e^{i\omega_l\zeta}}{\zeta}&\frac {1-e^{-i\zeta}}{\zeta}&\dots&
\frac {1-e^{-i\omega_l\zeta}}{\zeta}\\
\frac {1-e^{i\zeta}}{\zeta}&\omega_1\frac {1-e^{i\omega_1\zeta}}{\zeta}&\dots&
\omega_l\frac {1-e^{i\omega_l\zeta}}{\zeta}&(-1)\frac {1-e^{-i\zeta}}{\zeta}&
\dots&-\omega_l\frac {1-e^{-i\omega_l\zeta}}{\zeta}\\
\hdotsfor[3]{7} \\
\frac {1-e^{i\zeta}}{\zeta}&\omega_1^{2l}\frac {1-e^{i\omega_1\zeta}}{\zeta}&
\dots&\omega_l^{2l}\frac {1-e^{i\omega_l\zeta}}{\zeta}&
(-1)^{2l}\frac {1-e^{-i\zeta}}{\zeta}&\dots&
(-\omega_l)^{2l}\frac {1-e^{-i\omega_l\zeta}}{\zeta}\\
\end{bmatrix},
$$
we have
$$\left|\frac {{\mathfrak F}(\zeta)}{\zeta^{2l+1}} \right|_{\zeta=0}=
|{\mathfrak V}(1,\omega_1,\dots,\omega_{2l+1})|=(2l+2)^{l+1}.
$$
Since $|1-\omega_1|=2\sin\frac {\pi}{2l+2}$, this gives (3.9) after some
simplification.\hfill$\square$\medskip

{\bf Remark 8}. For $l=0$ (3.9) gives the classical formula for small ball
asymptotics of Brownian bridge under $L_2$-norm. For $l=1$ the formula (3.9)
was given in \cite[\S6]{BNO}, but the distortion constant there was
calculated only numerically.\medskip

Now we consider $m$-times integrated process
$(B_{\{l\}})_m^{[\beta_1,\,...,\,\beta_m]}(t)$. Due to
\cite[Theorem 2.1]{NN}, its covariance is the Green function of the BVP
$$\left\{
\begin{aligned}
&{\cal L}_{(B_{\{l\}})_m}u\equiv  (-1)^{l+m+1}u^{(2m+2l+2)}=\mu u\quad
\text{on} \quad [0,1],\\
&u(\beta_m)=u'(\beta_{m-1})=\dots=u^{(m-1)}(\beta_1)=0,\\
&u^{(m)}(0)=u^{(m)}(1)=0,\quad u^{(m+j)}(0)-u^{(m+j)}(1)=0,\ \ j=1,\dots,2l,\\
&u^{(m+2l+2)}(1-\beta_1)=u^{(m+2l+3)}(1-\beta_2)=\dots
=u^{(2m+2l+1)}(1-\beta_m)=0.\\
\end{aligned}
\right.\eqno(3.11)$$
The problem (3.11) satisfies all assumptions of Theorem 1.2 (with
$\ell=m+l+1$). This gives us the small ball asymptotics for the processes
$(B_{\{l\}})_m^{[\beta_1,\,...,\,\beta_m]}(t)$ up to a constant (note that
$\th_{\ell}=1$ and $\varkappa=(2m+2l+1)(m+l+1)-(2l+1)$). As for the
distortion constant, the only problem for its calculation is the length of
explicit representation of the corresponding Fredholm determinant. We
restrict ourselves to the case $l=1$.\medskip

{\bf Theorem 3.4}. { \it Let $m\in{\mathbb N}$. Then, as $\ep\to0$,
\begin{multline*}
{\bf P}\{\|(B_{\{1\}})_m^{[\beta_1,\,...,\,\beta_m]}\|\le\ep\}\sim\\
\sim\frac
{(2m+4)^{\frac {m+2}2}\sqrt {\textstyle2\sin\frac {3\pi\vphantom2}{2m+4}}}
{|{\mathfrak V}(z_{m+2}^{k_1},z_{m+2}^{k_2},
\dots,z_{m+2}^{k_m})|\cdot \sqrt{\prod \limits_{j=1}^m |1+z_{m+2}^{k_j}|^2+
\prod \limits_{j=1}^m |1+z_{m+2}^{k'_j}|^2}}\cdot\\
\cdot\frac {\ep_{m+2}^{-2}}{\sqrt {\pi {\mathfrak D}_{m+2}}}
\exp\left(-\ \frac {{\mathfrak D}_{m+2}}{2\ep_{m+2}^2}\right).
\tag{3.12}
\end{multline*}
where ${\mathfrak D}_{\ell}$ is defined in (2.4), $\ep_{\ell}$ is introduced
in Theorem 3.2, and for $j=1,\dots,m$}
$$k_j=\ \left\{
\begin{array}
{lll} m-j,& \mbox {if}& \beta_j=0,\\
m+3+j,& \mbox {if}&\beta_j=1,
\end{array}\right.
\qquad k'_j=2m+3-k_j.
\eqno(3.13)$$

\noindent {\bf Proof}. Put  $\zeta=\mu^{\frac 1{2m+4}}$. Substituting the
general solution of the equation (3.11) into boundary conditions, we deduce
that
$$C_{\rm dist}((B_{\{1\}})_m)=\prod_{n=1}^{\infty}\biggl(\frac {r_n}
{\pi (n-\frac {m-1}{2m+4})}\biggr)^{m+2},
$$
where $r_1<r_2<\ldots$ are positive roots of the entire function
\begin{multline*}
{\mathfrak F}_1(\zeta)\equiv\\
\det
\begin{bmatrix}
1&\omega_1^{k_1}&\dots&\omega_{m+1}^{k_1}&(-1)^{k_1}&\dots&
(-\omega_{m+1})^{k_1}\\
\hdotsfor[3]{7} \\
1&\omega_1^{k_m}&\dots&\omega_{m+1}^{k_m}&(-1)^{k_m}&\dots&
(-\omega_{m+1})^{k_m}\\
1&\omega_1^m&\dots&\omega_{m+1}^m&(-1)^m&\dots&(-\omega_{m+1})^m\\
1-e^{i\zeta}&\omega_1^{m+1}(1-e^{i\omega_1\zeta})&\dots&
\dots&(-1)^{m+1}(1-e^{-i\zeta})&\dots&\dots\\
1-e^{i\zeta}&\omega_1^{m+2}(1-e^{i\omega_1\zeta})&\dots&
\dots&(-1)^{m+2}(1-e^{-i\zeta})&\dots&\dots\\
e^{i\zeta}&\omega_1^me^{i\omega_1\zeta}&\dots&
\omega_{m+1}^me^{i\omega_{m+1}\zeta}&(-1)^me^{-i\zeta}&\dots&\dots\\
e^{i\zeta}&\omega_1^{k'_1}e^{i\omega_1\zeta}&\dots&
\omega_{m+1}^{k'_1}e^{i\omega_{m+1}\zeta}&(-1)^{k'_1}e^{-i\zeta}&\dots&\dots\\
\hdotsfor[3]{7} \\
e^{i\zeta}&\omega_1^{k'_m}e^{i\omega_1\zeta}&\dots&
\omega_{m+1}^{k'_m}e^{i\omega_{m+1}\zeta}&(-1)^{k'_m}e^{-i\zeta}&\dots&\dots\\
\end{bmatrix},
\end{multline*}
while $\omega_j=z_{m+2}^j$.

Similarly to Theorem 2.2 we obtain for $|\zeta|\to\infty$ and
$|\arg(\zeta)|\le \frac {\pi}{2m+4}$
$${\mathfrak F}_1(\zeta)=\exp(-i\omega_1\zeta)\exp(-i\omega_2\zeta)\dots
\exp(-i\omega_{m+1}\zeta)\cdot\Bigl(\Phi(\zeta)+O(|\zeta|^{-1})\Bigr),
$$

\noindent where
$$\Phi(\zeta)=\det
\begin{bmatrix}
1&\omega_1^{k_1}&\dots&\omega_{m+1}^{k_1}&(-1)^{k_1}&0&\dots&0\\
\hdotsfor[3]{8} \\
1&\omega_1^{k_m}&\dots&\omega_{m+1}^{k_m}&(-1)^{k_m}&0&\dots&0\\
1&\omega_1^m&\dots&\omega_{m+1}^m&(-1)^m&0&\dots&0\\
1-e^{i\zeta}&\omega_1^{m+1}&\dots&\omega_{m+1}^{m+1}&(-1)^{m+1}(1-e^{-i\zeta})
&(-\omega_1)^{m+1}&\dots&(-\omega_{m+1})^{m+1}\\
1-e^{i\zeta}&\omega_1^{m+2}&\dots&\omega_{m+1}^{m+2}&(-1)^{m+2}(1-e^{-i\zeta})
&(-\omega_1)^{m+2}&\dots&(-\omega_{m+1})^{m+2}\\
-e^{i\zeta}&0&\dots&0&(-1)^m(-e^{-i\zeta})&(-\omega_1)^m&\dots&
(-\omega_{m+1})^m\\
-e^{i\zeta}&0&\dots&0&(-1)^{k'_1}(-e^{-i\zeta})&(-\omega_1)^{k'_1}&\dots&
(-\omega_{m+1})^{k'_1}\\
\hdotsfor[3]{8} \\
-e^{i\zeta}&0&\dots&0&(-1)^{k'_m}(-e^{-i\zeta})&(-\omega_1)^{k'_m}&\dots&
(-\omega_{m+1})^{k'_m}\\
\end{bmatrix}.$$

Expanding this determinant in the elements of the first and $(m+3)$-rd
columns we derive
$$|\Phi(\zeta)|={\mathfrak M}\cdot
|\exp(-i\zeta)+\omega_1^{-3}\exp(i\zeta)+R|,$$
where
\begin{multline*}
{\mathfrak M}=
|{\mathfrak V}(\omega_1^{k_1},\dots,\omega_1^{k_m},\omega_1^m,\omega_1^{m+1})
\cdot {\mathfrak V}
(\omega_1^{m+2},\omega_1^m,\omega_1^{k'_1},\dots,\omega_1^{k'_m})+\\
+{\mathfrak V}(\omega_1^{k_1},\dots,\omega_1^{k_m},\omega_1^m,\omega_1^{m+2})
\cdot {\mathfrak V}
(\omega_1^{m+1},\omega_1^m,\omega_1^{k'_1},\dots,\omega_1^{k'_m})|,
\end{multline*}
while $R$ is a constant which is inessential for us.

Setting $\delta=-\frac {m-1}{2m+4}$ we obtain in view of
$|{\mathfrak F}_1(\zeta)|\equiv |{\mathfrak F}_1(\omega_1\zeta)|$
$$\frac {|{\mathfrak F}_1(\zeta)|}
{\Bigl|\zeta^3\prod \limits_{j=0}^{m+1}\psi_{\delta}(\omega_j\zeta)\Bigr|}\to
\frac {2^{m+2}\pi^{m-1}{\mathfrak M}}{\Gamma^{2m+4}(1+\delta)},
$$
for $\zeta\to\infty$, $\arg(\zeta)\ne \frac{\pi j}{2m+4}$, $j\in\mathbb Z$.

This implies, similarly to the proof of Theorem 2.2,
$$C^2_{\rm dist}((B_{\{1\}})_m)=
\frac {\Gamma^{2m+4}(1+\delta)}
{2^{m+2}\pi^{m-1}{\mathfrak M}}\cdot
\left|\frac {{\mathfrak F}_1(\zeta)}{\zeta^3} \right|_{\zeta=0}.
$$
Subtracting in the determinant $(m+1)$-st row from $(m+4)$-th one we obtain
in view of (3.13)
$$\left|\frac {{\mathfrak F}_1(\zeta)}{\zeta^3} \right|_{\zeta=0}=
|{\mathfrak V}(1,\omega_1,\omega_1^2,\dots,\omega_1^{2m+3})|=(2m+4)^{m+2}.
$$
It remains to take into account that due to (3.13)
$${\mathfrak M}=
|{\mathfrak V}(\omega_1^{k_1},\dots,\omega_1^{k_m})|^2\cdot
\frac {2m+4}{|1-\omega_3|}\cdot
\Big(\prod \limits_{j=1}^m |1+\omega_1^{k_j}|^2+
\prod \limits_{j=1}^m |1+\omega_1^{k'_j}|^2\Big).
$$
Since $|1-\omega_3|=2\sin\frac {3\pi}{2m+4}$, we arrive at (3.12) after some
simplification.\hfill$\square$\medskip

\section{Integrated centered Wiener process and \\ related processes}

Similarly to Section 3, we define the sequence of integrated centered
analogues of Wiener process. We set
$$W_{\{0\}}(t)=W(t); \qquad
W_{\{l\}}(t)=\int_0^t\overline {W_{\{l-1\}}}(s)\,ds, \quad
l\in\mathbb N.$$
The spectrum of the process $W_{\{1\}}$ and its $L_2$-small ball asymptotics
was studied in \cite[\S7]{BNO}. In \cite[Example 5.4]{NN} it is pointed out
that the covariance $G_{W_{\{1\}}}$ is the Green function of the BVP
$${\cal L}_{W_{\{1\}}}u\equiv u^{\scriptscriptstyle IV}=\mu u\quad
\mbox{\rm on} \quad [0,1], \qquad u(0)=u(1)=u''(0)=u''(1)=0.\eqno(4.1)$$
Theorem 3.1 allows us to write down the BVPs generated by processes
$W_{\{l\}}$ and $\overline {W_{\{l\}}}$. We begin from the second process
and prove an unexpected relation.\medskip

{\bf Theorem 4.1}. { \it Let $l\in{\mathbb N}_0$. Then the following
distributional equality holds true:}
$$\|\overline {W_{\{l\}}}\|\stackrel {d}{=}\|B_{\{l\}}\|.\eqno(4.2)$$

\noindent {\bf Proof}. For $l=0$ the equality (4.2) is well known, see, e.g.,
\cite{DMY} and \cite[\S3]{BNO}. Let $l\ge1$. Applying in turns the second and
the first statements of Theorem 3.1 to the problem (4.1) we deduce that the
covariance $G_{\overline {W_{\{l\}}}}(t,s)$ is the Green function of the BVP
$$\left\{
\begin{aligned}
&{\cal L}_{\overline {W_{\{l\}}}}u\equiv (-1)^{l+1}u^{(2l+2)}=\mu u\quad
\text{on} \quad [0,1],\\
&u^{(l+1)}(0)=u^{(l+1)}(1)=0,\quad u^{(j)}(0)-u^{(j)}(1)=0,\ \
j=0,\dots,l-1,\,l+2,\dots,2l+1.\\
\end{aligned}
\right.\eqno(4.3)$$
It is easy to check that $(l+1)$-times differentiation maps mutually the
eigenfunctions of BVPs (3.10) and (4.3) (if the corresponding eigenvalue
$\mu\ne0$). Hence nonzero eigenvalues of these BVPs coincide pairwise, and
therefore nonzero eigenvalues of the covariances also coincide. This gives
(4.2).\hfill$\square$\medskip

{\bf Theorem 4.2}. {\it Let $l\in{\mathbb N}$. Then, as $\ep\to0$,

$${\bf P}\{\|W_{\{l\}}\|\le\ep\}\sim
(2l+2)^{\frac 32}\,{\textstyle\sin\frac {\pi}{2l+2}\,
\langle\cos\frac {\pi}{2l+2}\rangle}\cdot
\frac {\ep_{l+1}^{-(2l-1)}}{\sqrt {\pi {\mathfrak D}_{l+1}}}
\exp\left(-\ \frac {{\mathfrak D}_{l+1}}{2\ep_{l+1}^2}\right)
\eqno(4.4)$$
with the same notation as in Theorem 3.2 (the angle brackets must be omitted
if $l$ is even).}\medskip

{\bf Remark 9}. For $l=1$ (4.4) coincides with formula obtained in
\cite[\S6]{BNO}, see also \cite[Proposition 1.7]{Na}.\medskip

\noindent {\bf Proof}. Similarly to Theorem 4.1, the covariance
$G_{W_{\{l\}}}(t,s)$ is the Green function of the BVP
$$\left\{
\begin{aligned}
&{\cal L}_{W_{\{l\}}}u\equiv (-1)^{l+1}u^{(2l+2)}=\mu u\quad
\text{on} \quad [0,1],\\
&u(0)=u(1)=0,\quad u^{(l+1)}(0)=u^{(l+1)}(1)=0,\\
&u^{(j)}(0)-u^{(j)}(1)=0,\quad j=1,\dots,l-1,\,l+2,\dots,2l.\\
\end{aligned}
\right.\eqno(4.5)$$
Since the problem (4.5) satisfies all assumptions of Theorem 1.2 (with
$\ell=l+1$), to prove (4.4) we only need to calculate the distortion
constant. Note that $\th_{\ell}=1$ and $\varkappa=(2l+1)l+1$.

Put $\zeta=\mu^{\frac 1{2l+2}}$. Substituting the general solution of the
equation (4.5) into boundary conditions, we deduce that
$$C_{\rm dist}(W_{\{l\}})=\prod_{n=1}^{\infty}\biggl(\frac {r_n}
{\pi (n+\frac {l-1}{2l+2})}\biggr)^{l+1},
$$
where $r_1<r_2<\ldots$ are positive roots of the entire function
$${\mathbb F}(\zeta)\equiv\det
\begin{bmatrix}
1&1&\dots&1&1&\dots&1\\
e^{i\zeta}&e^{i\omega_1\zeta}&\dots&
e^{i\omega_l\zeta}&e^{-i\zeta}&\dots&e^{-i\omega_l\zeta}\\
1-e^{i\zeta}&\omega_1(1-e^{i\omega_1\zeta})&\dots&
\dots&(-1)(1-e^{-i\zeta})&\dots&\dots\\
1-e^{i\zeta}&\omega_1^2(1-e^{i\omega_1\zeta})&\dots&
\dots&(-1)^2(1-e^{-i\zeta})&\dots&\dots\\
\hdotsfor[3]{7} \\
1-e^{i\zeta}&\omega_1^{l-1}(1-e^{i\omega_1\zeta})&\dots&
\dots&(-1)^{l-1}(1-e^{-i\zeta})&\dots&\dots\\
1&-1&\dots&(-1)^l&(-1)^{l+1}&\dots&-1\\
e^{i\zeta}&-e^{-i\omega_1\zeta}&\dots&
(-1)^le^{-i\omega_l\zeta}&(-1)^{l+1}e^{-i\zeta}&\dots&-e^{-i\omega_l\zeta}\\
1-e^{i\zeta}&\omega_1^{l+2}(1-e^{i\omega_1\zeta})&\dots&
\dots&(-1)^{l+2}(1-e^{-i\zeta})&\dots&\dots\\
\hdotsfor[3]{7} \\
1-e^{i\zeta}&\omega_1^{2l}(1-e^{i\omega_1\zeta})&\dots&
\dots&(-1)^{2l}(1-e^{-i\zeta})&\dots&\dots\\
\end{bmatrix},
$$
while $\omega_j=z_{l+1}^j$.

Subtracting the first row from the second one and $(l+2)$-nd row from
$(l+3)$-rd one, similarly to Theorem 2.2 we obtain for $|\zeta|\to\infty$ and
$|\arg(\zeta)|\le \frac {\pi}{2l+2}$
$${\mathbb F}(\zeta)=\exp(-i\omega_1\zeta)\exp(-i\omega_2\zeta)
\dots\exp(-i\omega_l\zeta)\cdot\Bigl(\Phi(\zeta)+O(|\zeta|^{-1})\Bigr),
$$

\noindent where
\begin{multline*}
\Phi(\zeta)=\det\\
\begin{bmatrix}
1&1&\dots&1&1&0&\dots&0\\
1-e^{i\zeta}&1&\dots&1&1-e^{-i\zeta}
&-1&\dots&-1\\
1-e^{i\zeta}&\omega_1&\dots&\omega_l&(-1)(1-e^{-i\zeta})
&\omega_1&\dots&\omega_l\\
\hdotsfor[3]{8} \\
1-e^{i\zeta}&\omega_1^{l-1}&\dots&\omega_l^{l-1}&(-1)^{l-1}(1-e^{-i\zeta})
&-(-\omega_1)^{l-1}&\dots&-(-\omega_l)^{l-1}\\
1&-1&\dots&(-1)^l&(-1)^{l+1}&0&\dots&0\\
1-e^{i\zeta}&-1&\dots&(-1)^l&(-1)^{l+1}(1-e^{-i\zeta})
&-(-1)^{l+2}&\dots&1\\
1-e^{i\zeta}&-\omega_1&\dots&(-1)^l\omega_l&(-1)^{l+2}(1-e^{-i\zeta})
&(-1)^{l+2}\omega_1&\dots&-\omega_l\\
\hdotsfor[3]{8} \\
1-e^{i\zeta}&-\omega_1^{l-1}&\dots&(-1)^l\omega_l^{l-1}&
(-1)^{2l}(1-e^{-i\zeta})&(-1)^{l+3}(-\omega_1)^{l-1}&\dots&(-\omega_l)^{l-1}\\
\end{bmatrix}.
\end{multline*}

\noindent Adding the upper half of the matrix to the lower one we get
$|\Phi(\zeta)|=2^{l+1}|\Delta_1(\zeta)|\cdot|\Delta_2(\zeta)|$, where

for even $l$
$$\Delta_1(\zeta)=\det
\begin{bmatrix}
1&1&\dots&1&0&\dots&0\\
1-e^{i\zeta}&1&\dots&1&1&\dots&1\\
1-e^{i\zeta}&\omega_2&\dots&\omega_l&\omega_{l+2}&\dots&\omega_{2l}\\
\hdotsfor[3]{7} \\
1-e^{i\zeta}&\omega_2^{l-1}&\dots&\omega_l^{l-1}&\omega_{l+2}^{l-1}&\dots&
\omega_{2l}^{l-1}\\
\end{bmatrix},
$$
$$\Delta_2(\zeta)=\det
\begin{bmatrix}
1&0&\dots&0&1&\dots&1\\
1-e^{-i\zeta}&1&\dots&1&1&\dots&1\\
1-e^{-i\zeta}&\omega_2&\dots&\omega_l&\omega_{l+2}&\dots&\omega_{2l}\\
\hdotsfor[3]{7} \\
1-e^{-i\zeta}&\omega_2^{l-1}&\dots&\omega_l^{l-1}&\omega_{l+2}^{l-1}&\dots&
\omega_{2l}^{l-1}\\
\end{bmatrix};
$$

for odd $l$
$$\Delta_1(\zeta)=\det
\begin{bmatrix}
1&1&\dots&1&1&0&\dots&0\\
1-e^{i\zeta}&1&\dots&1&1-e^{-i\zeta}
&1&\dots&1\\
1-e^{i\zeta}&\omega_2&\dots&\omega_{l-1}&(-1)(1-e^{-i\zeta})
&-\omega_2&\dots&-\omega_{l-1}\\
\hdotsfor[3]{8} \\
1-e^{i\zeta}&\omega_2^{l-1}&\dots&\omega_{l-1}^{l-1}&(-1)^{l-1}(1-e^{-i\zeta})
&(-\omega_2)^{l-1}&\dots&(-\omega_{l-1})^{l-1}\\
\end{bmatrix},
$$
$$\Delta_2(\zeta)=\det
\begin{bmatrix}
1&1&\dots&1&0&\dots&0\\
1&1&\dots&1&1&\dots&1\\
\omega_1&\omega_3&\dots&\omega_l&-\omega_1&\dots&-\omega_l\\
\hdotsfor[3]{7} \\
\omega_1^{l-1}&\omega_3^{l-1}&\dots&\omega_l^{l-1}&(-\omega_1)^{l-1}&\dots&
(-\omega_l)^{l-1}\\
\end{bmatrix}
$$
(note that in this case $\Delta_1$ coincides with the determinant from
Theorem 3.3).

Expanding these determinants in the elements of the first row we obtain:

for even $l$
\begin{multline*}
|\Phi(\zeta)|=\frac
{2^{l+1}|{\mathfrak V}(1,\omega_2,\dots,\omega_{2l-2})|^2}{|1-\omega_1|^2}
\cdot |1-\omega_1\exp(i\zeta)|^2=\\
=\frac {4(2l+2)^{l-1}}{|1-\omega_1|^2}\cdot
|\exp(-i\zeta)+\omega_1^2\exp(i\zeta)-2\omega_1|;
\end{multline*}

for odd $l$
\begin{multline*}
|\Phi(\zeta)|=\frac
{2^{l+3}|{\mathfrak V}(1,\omega_2,\dots,\omega_{2l-2})|^2}{|1-\omega_2|^2}
\cdot |(1-\exp(i\zeta))(\exp(-i\zeta)-\omega_1^2)|=\\
=\frac {16(2l+2)^{l-1}}{|1-\omega_2|^2}\cdot
|\exp(-i\zeta)+\omega_1^2\exp(i\zeta)-(1+\omega_1^2)|.
\end{multline*}

Setting $\delta=\frac {l-1}{2l+2}$ we obtain in view of
$|{\mathbb F}(\zeta)|\equiv |{\mathbb F}(\omega_1\zeta)|$
$$\frac {|{\mathbb F}(\zeta)|}
{|\zeta^{2l}\prod \limits_{j=0}^{l}\psi_{\delta}(\omega_j\zeta)|}\to
\frac {2^{l+3}(2l+2)^{l-1}}
{\Gamma^{2l+2}(1+\delta)\pi^{l-1}\,{\mathbb M}}
$$
for $\zeta\to\infty$, $\arg(\zeta)\ne \frac{\pi j}{2l+2}$, $j\in\mathbb Z$
(here ${\mathbb M}=|1-\omega_1|^2$ for even $l$ and
${\mathbb M}=|1-\omega_2|^2/4$ for odd $l$).

This implies, similarly to the proof of Theorem 2.2,
$$C^2_{\rm dist}(W_{\{l\}})=
\frac {\Gamma^{2l+2}(1+\delta)\pi^{l-1}\,{\mathbb M}}
{2^{l+3}(2l+2)^{l-1}}\cdot
\left|\frac {{\mathbb F}(\zeta)}{\zeta^{2l}} \right|_{\zeta=0}.
$$

Since
$$\frac {{\mathbb F}(\zeta)}{\zeta^{2l}}=
\det
\begin{bmatrix}
1&1&\dots&1&1&\dots&1\\
\frac {1-e^{i\zeta}}{\zeta}&\frac {1-e^{i\omega_1\zeta}}{\zeta}&\dots&
\frac {1-e^{i\omega_l\zeta}}{\zeta}&\frac {1-e^{-i\zeta}}{\zeta}&\dots&
\frac {1-e^{-i\omega_l\zeta}}{\zeta}\\
\frac {1-e^{i\zeta}}{\zeta}&\omega_1\frac {1-e^{i\omega_1\zeta}}{\zeta}&\dots&
\omega_l\frac {1-e^{i\omega_l\zeta}}{\zeta}&(-1)\frac {1-e^{-i\zeta}}{\zeta}&
\dots&-\omega_l\frac {1-e^{-i\omega_l\zeta}}{\zeta}\\
\hdotsfor[3]{7} \\
\frac {1-e^{i\zeta}}{\zeta}&\omega_1^{l-1}\frac {1-e^{i\omega_1\zeta}}{\zeta}&
\dots&\omega_l^{l-1}\frac {1-e^{i\omega_l\zeta}}{\zeta}&
(-1)^{l-1}\frac {1-e^{-i\zeta}}{\zeta}&\dots&
(-\omega_l)^{l-1}\frac {1-e^{-i\omega_l\zeta}}{\zeta}\\
1&-1&\dots&(-1)^{l\vphantom{\frac12}}&(-1)^{l+1}&\dots&-1\\
\frac {1-e^{i\zeta}}{\zeta}&-\frac {1-e^{i\omega_1\zeta}}{\zeta}&\dots&
(-1)^l\frac {1-e^{i\omega_l\zeta}}{\zeta}&
(-1)^{l+1}\frac {1-e^{-i\zeta}}{\zeta}&\dots&
-\frac {1-e^{-i\omega_l\zeta}}{\zeta}\\
\frac {1-e^{i\zeta}}{\zeta}&\omega_1^{l+2}\frac {1-e^{i\omega_1\zeta}}{\zeta}
&\dots&\omega_l^{l+2\vphantom{\frac12}}\frac {1-e^{i\omega_l\zeta}}{\zeta}&
(-1)^{l+2}\frac {1-e^{-i\zeta}}{\zeta}&
\dots&(-\omega_l)^{l+2}\frac {1-e^{-i\omega_l\zeta}}{\zeta}\\
\hdotsfor[3]{7} \\
\frac {1-e^{i\zeta}}{\zeta}&\omega_1^{2l}\frac {1-e^{i\omega_1\zeta}}{\zeta}&
\dots&\omega_l^{2l}\frac {1-e^{i\omega_l\zeta}}{\zeta}&
(-1)^{2l}\frac {1-e^{-i\zeta}}{\zeta}&\dots&
(-\omega_l)^{2l}\frac {1-e^{-i\omega_l\zeta}}{\zeta}\\
\end{bmatrix},
$$
we have
$$\left|\frac {{\mathbb F}(\zeta)}{\zeta^{2l}} \right|_{\zeta=0}=
|{\mathfrak V}(1,\omega_1,\dots,\omega_{2l+1})|=(2l+2)^{l+1},
$$
Since $|1-\omega_j|=2\sin\frac {j\pi}{2l+2}$, this gives (4.4) after some
simplification.\hfill$\square$\medskip

{\bf Remark 9}. In our proof we in fact use that the operator of the BVP
(4.5) is the square of the BVP operator of order $l+1$
$$\left\{
\begin{aligned}
&{\mathfrak L}u\equiv i^{l+1}u^{(l+1)}=\mu u\quad \text{on} \quad [0,1],\\
&u(0)=u(1)=0,\quad u^{(j)}(0)-u^{(j)}(1)=0,\quad j=1,\dots,l-1.\\
\end{aligned}
\right.$$
For odd $l$ the operator ${\mathfrak L}$ coincides with the operator
${\cal L}_{B_{\{\frac {l-1}2\}}}$, see (3.10).\medskip

Now we consider $m$-times integrated process
$(W_{\{l\}})_m^{[\beta_1,\,...,\,\beta_m]}(t)$. According to
\cite[Theorem 2.1]{NN}, its covariance is the Green function of the BVP
$$\left\{
\begin{aligned}
&{\cal L}_{(W_{\{l\}})_m}u\equiv  (-1)^{l+m+1}u^{(2m+2l+2)}=\mu u\quad
\text{on} \quad [0,1],\\
&u(\beta_m)=u'(\beta_{m-1})=\dots=u^{(m-1)}(\beta_1)=0,\\
&u^{(m)}(0)=u^{(m)}(1)=0,\qquad u^{(m+j)}(0)-u^{(m+j)}(1)=0,\ \ j=1,\dots,l-1,\\
&u^{(m+l+1)}(0)=u^{(m+l+1)}(1)=0,\quad u^{(m+j)}(0)-u^{(m+j)}(1)=0,\ \
j=l+2,\dots,2l,\\
&u^{(m+2l+2)}(1-\beta_1)=u^{(m+2l+3)}(1-\beta_2)=\dots
=u^{(2m+2l+1)}(1-\beta_m)=0.\\
\end{aligned}
\right.\eqno(4.6)$$
The problem (4.6) satisfies all assumptions of Theorem 1.2 (with
$\ell=m+l+1$). This gives us the small ball asymptotics for the processes
$(W_{\{l\}})_m^{[\beta_1,\,...,\,\beta_m]}(t)$ up to constant (note that
$\th_{\ell}=1$ and $\varkappa=(2m+2l+1)(m+l+1)-2l$). For $l=1$ the problem
(4.6) has separated boundary conditions; distortion constants in this case
were calculated in \cite[Proposition 1.7]{Na}. We restrict ourselves to the
case $l=2$.\medskip

{\bf Theorem 4.3}. { \it Let $m\in{\mathbb N}$. Then, as $\ep\to0$,
\begin{multline*}
{\bf P}\{\|(W_{\{2\}})_m^{[\beta_1,\,...,\,\beta_m]}\|\le\ep\}\sim\\
\sim\frac
{4(2m+6)^{\frac {m+2}2}\,{\textstyle\sin\frac {\pi}{m+3}}\,
\sqrt {\textstyle\sin\frac {\pi\vphantom2}{2m+6}\sin\frac {5\pi}{2m+6}}}
{|{\mathfrak V}(z_{m+3}^{k_1},z_{m+3}^{k_2},
\dots,z_{m+3}^{k_m})|\cdot \sqrt{\prod \limits_{j=1}^m
|z_{m+3}+z_{m+3}^{k_j}|^2+\prod \limits_{j=1}^m
|z_{m+3}+z_{m+3}^{k'_j}|^2}}\cdot\\
\cdot\frac {\ep_{m+3}^{-3}}{\sqrt {\pi {\mathfrak D}_{m+3}}}
\exp\left(-\ \frac {{\mathfrak D}_{m+3}}{2\ep_{m+3}^2}\right),
\tag{4.7}
\end{multline*}
where ${\mathfrak D}_{\ell}$ is defined in (2.4), $\ep_{\ell}$ is introduced
in Theorem 3.2, and for $j=1,\dots,m$}
$$k_j=\ \left\{
\begin{array}
{lll} m-j,& \mbox {if}& \beta_j=0,\\
m+5+j,& \mbox {if}&\beta_j=1,
\end{array}\right.
\qquad k'_j=2m+5-k_j.
\eqno(4.8)$$

\noindent {\bf Proof}. Put $\zeta=\mu^{\frac 1{2m+6}}$. Substituting the
general solution of the equation (4.6) into boundary conditions, we deduce
that
$$C_{\rm dist}((W_{\{2\}})_m)=\prod_{n=1}^{\infty}\biggl(\frac {r_n}
{\pi (n-\frac {m-1}{2m+6})}\biggr)^{m+3},
$$
where $r_1<r_2<\ldots$ are positive roots of the entire function
\begin{multline*}
{\mathbb F}_1(\zeta)\equiv\\
\det
\begin{bmatrix}
1&\omega_1^{k_1}&\dots&\omega_{m+2}^{k_1}&(-1)^{k_1}&\dots&
(-\omega_{m+2})^{k_1}\\
\hdotsfor[3]{7} \\
1&\omega_1^{k_m}&\dots&\omega_{m+2}^{k_m}&(-1)^{k_m}&\dots&
(-\omega_{m+2})^{k_m}\\
1&\omega_1^m&\dots&\omega_{m+2}^m&(-1)^m&\dots&(-\omega_{m+2})^m\\
1&\omega_1^{m+3}&\dots&\omega_{m+2}^{m+3\vphantom{\frac12}}&(-1)^{m+3}&\dots&
(-\omega_{m+2})^{m+3}\\
1-e^{i\zeta}&\omega_1^{m+1}(1-e^{i\omega_1\zeta})&\dots&
\dots&(-1)^{m+1}(1-e^{-i\zeta})&\dots&\dots\\
1-e^{i\zeta}&\omega_1^{m+4}(1-e^{i\omega_1\zeta})&\dots&
\dots&(-1)^{m+4}(1-e^{-i\zeta})&\dots&\dots\\
e^{i\zeta}&\omega_1^me^{i\omega_1\zeta}&\dots&
\omega_{m+2}^me^{i\omega_{m+2}\zeta}&(-1)^me^{-i\zeta}&\dots&\dots\\
e^{i\zeta}&\omega_1^{m+3}e^{i\omega_1\zeta}&\dots&
\omega_{m+2}^{m+3\vphantom{\frac12}}e^{i\omega_{m+2}\zeta}&
(-1)^{m+3}e^{-i\zeta}&\dots&\dots\\
e^{i\zeta}&\omega_1^{k'_1}e^{i\omega_1\zeta}&\dots&
\omega_{m+2}^{k'_1\vphantom{\frac12}}e^{i\omega_{m+1}\zeta}&
(-1)^{k'_1}e^{-i\zeta}&\dots&\dots\\
\hdotsfor[3]{7} \\
e^{i\zeta}&\omega_1^{k'_m}e^{i\omega_1\zeta}&\dots&
\omega_{m+2}^{k'_m}e^{i\omega_{m+1}\zeta}&(-1)^{k'_m}e^{-i\zeta}&\dots&\dots\\
\end{bmatrix},
\end{multline*}
while $\omega_j=z_{m+3}^j$.

Similarly to Theorem 2.2 we obtain for $|\zeta|\to\infty$ and
$|\arg(\zeta)|\le \frac {\pi}{2m+6}$
$${\mathbb F}_1(\zeta)=\exp(-i\omega_1\zeta)\exp(-i\omega_2\zeta)\dots
\exp(-i\omega_{m+2}\zeta)\cdot\Bigl(\Phi(\zeta)+O(|\zeta|^{-1})\Bigr),
$$

\noindent where
$$\Phi(\zeta)=\det
\begin{bmatrix}
1&\omega_1^{k_1}&\dots&\omega_{m+1}^{k_1}&(-1)^{k_1}&0&\dots&0\\
\hdotsfor[3]{8} \\
1&\omega_1^{k_m}&\dots&\omega_{m+2}^{k_m}&(-1)^{k_m}&0&\dots&0\\
1&\omega_1^m&\dots&\omega_{m+2}^m&(-1)^m&0&\dots&0\\
1&\omega_1^{m+3}&\dots&\omega_{m+2}^{m+3}&(-1)^{m+3}&0&\dots&0\\
1-e^{i\zeta}&\omega_1^{m+1}&\dots&\omega_{m+2}^{m+1}&(-1)^{m+1}(1-e^{-i\zeta})
&(-\omega_1)^{m+1}&\dots&(-\omega_{m+2})^{m+1}\\
1-e^{i\zeta}&\omega_1^{m+4}&\dots&\omega_{m+2}^{m+4}&(-1)^{m+4}(1-e^{-i\zeta})
&(-\omega_1)^{m+4}&\dots&(-\omega_{m+2})^{m+4}\\
-e^{i\zeta}&0&\dots&0&(-1)^m(-e^{-i\zeta})&(-\omega_1)^m&\dots&
(-\omega_{m+2})^m\\
-e^{i\zeta}&0&\dots&0&(-1)^{m+3}(-e^{-i\zeta})&(-\omega_1)^{m+3}&\dots&
(-\omega_{m+2})^{m+3}\\
-e^{i\zeta}&0&\dots&0&(-1)^{k'_1}(-e^{-i\zeta})&(-\omega_1)^{k'_1}&\dots&
(-\omega_{m+2})^{k'_1}\\
\hdotsfor[3]{8} \\
-e^{i\zeta}&0&\dots&0&(-1)^{k'_m}(-e^{-i\zeta})&(-\omega_1)^{k'_m}&\dots&
(-\omega_{m+2})^{k'_m}\\
\end{bmatrix}.$$

Expanding this determinant in the elements of the first and $(m+4)$-rd
columns we derive
$$|\Phi(\zeta)|={\mathbb M}_1\cdot
|\exp(-i\zeta)+\omega_1^{-4}\exp(i\zeta)+R|,$$
where
\begin{multline*}
{\mathbb M}_1=|{\mathfrak V}(\omega_1^{k_1},\dots,
\omega_1^{k_m},\omega_1^m,\omega_1^{m+3},\omega_1^{m+1})\cdot
{\mathfrak V}(\omega_1^{m+4},\omega_1^m,\omega_1^{m+3},\omega_1^{k'_1},
\dots,\omega_1^{k'_m})+\\
+{\mathfrak V}(\omega_1^{k_1},\dots,
\omega_1^{k_m},\omega_1^m,\omega_1^{m+3},\omega_1^{m+4})\cdot
{\mathfrak V}(\omega_1^{m+1},\omega_1^m,\omega_1^{m+3},\omega_1^{k'_1},
\dots,\omega_1^{k'_m})|,
\end{multline*}
while $R$ is a constant which is inessential for us.

Setting $\delta=-\frac {m-1}{2m+6}$ we obtain in view of
$|{\mathbb F}_1(\zeta)|\equiv |{\mathbb F}_1(\omega_1\zeta)|$
$$\frac {|{\mathbb F}_1(\zeta)|}
{\Bigl|\zeta^4\prod \limits_{j=0}^{m+2}\psi_{\delta}(\omega_j\zeta)\Bigr|}\to
\frac {2^{m+3}\pi^{m-1}{\mathbb M}_1}{\Gamma^{2m+6}(1+\delta)},
$$
for $\zeta\to\infty$, $\arg(\zeta)\ne \frac{\pi j}{2m+6}$, $j\in\mathbb Z$.

This implies, similarly to the proof of Theorem 2.2,
$$C^2_{\rm dist}((W_{\{2\}})_m)=
\frac {\Gamma^{2m+6}(1+\delta)}
{2^{m+3}\pi^{m-1}{\mathbb M}_1}\cdot
\left|\frac {{\mathbb F}_1(\zeta)}{\zeta^4} \right|_{\zeta=0}.
$$
Let us subtract in the determinant $(m+1)$-st row from $(m+5)$-th one and
$(m+2)$-nd row from $(m+6)$-th one. In view of (4.8) we obtain
$$\left|\frac {{\mathbb F}_1(\zeta)}{\zeta^4} \right|_{\zeta=0}=
|{\mathfrak V}(1,\omega_1,\omega_1^2,\dots,\omega_1^{2m+5})|=(2m+6)^{m+3}.
$$
It remains to take into account that due to (4.8)
$${\mathbb M}_1=
|{\mathfrak V}(\omega_1^{k_1},\dots,\omega_1^{k_m})|^2\cdot
\frac {(2m+6)^2}{|1-\omega_1||1-\omega_2|^2|1-\omega_5|}\cdot
\Big(\prod \limits_{j=1}^m |\omega_1+\omega_1^{k_j}|^2+
\prod \limits_{j=1}^m |\omega_1+\omega_1^{k'_j}|^2\Big).
$$
Since $|1-\omega_j|=2\sin\frac {j\pi}{2m+6}$ we arrive at (4.7) after some
simplification.\hfill$\square$\medskip

\bigskip
I am grateful to Professor Ya.Yu.~Nikitin for some comments and references,
and also for his constant encouragement.

\end{document}